\documentclass[a4paper,12pt]{amsart}
\usepackage{amssymb}
\usepackage{ifthen}
\usepackage[usenames]{color}
\usepackage{graphicx}
\nonstopmode \numberwithin{equation}{section}
\newtheorem{thm}{Theorem}[section]

\newtheorem{cor}{Corollary}[section]
\newtheorem{lem}{Lemma}[section]
\newtheorem{prop}{Proposition}[section]

\newtheorem{claim}{Claim}[section]
\newtheorem{subclaim}{Subclaim}
\newtheorem{conj}[equation]{Conjecture}
\newtheorem{case}{Case}[section]
\newtheorem*{mysolution}{Solution}
\newtheorem{step}{Step}[section]
\theoremstyle{definition}
\newtheorem{defn}{Definition}[section]
\newtheorem{examp}{Example}[section]
\newtheorem{prob}[equation]{Problem}
\newtheorem{ques}[equation]{Question}
\newtheorem{rem}{Remark}[section]
\newcounter {own}
\def\theown {\thesection       .\arabic{own}}

\newenvironment{pf}[1][]{%
 \vskip 3mm
 \noindent
 \ifthenelse{\equal{#1}{}}%
  {{\slshape Proof. }}%
  {{\slshape #1.} }%
 }%
{\qed\bigskip}

\newcounter{alphabet}
\newcounter{tmp}
\newenvironment{Thm}[1][]{\refstepcounter{alphabet}%
\bigskip%
\noindent%
{\bf Theorem \Alph{alphabet}}%
\ifthenelse{\equal{#1}{}}{}{ (#1)}%
{\bf .} \itshape}{\vskip 8pt}

\makeatletter

\renewcommand{\Ref}[1]{\@ifundefined{r@#1}{}{\setcounter{tmp}{\ref{#1}}\Alph{tmp}}}

\newcommand{\IR}{{\mathbb R}}

\newcommand{\IC}{{\mathbb C}}
\newcommand{\IS}{{\mathbb S}}
\newcommand{\ID}{{\mathbb D}}
\newcommand{\IB}{{\mathbb B}}

\def\be{\begin{equation}}
\def\ee{\end{equation}}

\newcommand{\ben}{\begin{enumerate}}
\newcommand{\een}{\end{enumerate}}

\newcommand{\blem}{\begin{lem}}
\newcommand{\elem}{\end{lem}}
\newcommand{\bthm}{\begin{thm}}
\newcommand{\ethm}{\end{thm}}
\newcommand{\bcor}{\begin{cor}}
\newcommand{\ecor}{\end{cor}}
\newcommand{\beg}{\begin{examp}}
\newcommand{\eeg}{\end{examp}}
\newcommand{\begs}{\begin{examples}}
\newcommand{\eegs}{\end{examples}}
\newcommand{\bdefe}{\begin{defn}}
\newcommand{\edefe}{\end{defn}}
\newcommand{\bprob}{\begin{prob}}
\newcommand{\eprob}{\end{prob}}
\newcommand{\bques}{\begin{ques}}
\newcommand{\eques}{\end{ques}}
\newcommand{\bei}{\begin{itemize}}
\newcommand{\eei}{\end{itemize}}
\newcommand{\bcl}{\begin{claim}}
\newcommand{\ecl}{\end{claim}}
\newcommand{\bscl}{\begin{subclaim}}
\newcommand{\escl}{\end{subclaim}}
\newcommand{\bca}{\begin{case}}
\newcommand{\eca}{\end{case}}
\newcommand{\bstep}{\begin{step}}
\newcommand{\estep}{\end{step}}
\newcommand{\bsol}{\begin{mysolution}}
\newcommand{\esol}{\end{mysolution}}
\newcommand{\bcon}{\begin{conj}}
\newcommand{\econ}{\end{conj}}
\newcommand{\bcons}{\begin{conjs}}
\newcommand{\econs}{\end{conjs}}
\newcommand{\bprop}{\begin{prop}}
\newcommand{\eprop}{\end{prop}}
\newcommand{\br}{\begin{rem}}
\newcommand{\er}{\end{rem}}
\newcommand{\brs}{\begin{rems}}
\newcommand{\ers}{\end{rems}}
\newcommand{\bo}{\begin{obser}}
\newcommand{\eo}{\end{obser}}
\newcommand{\bos}{\begin{obsers}}
\newcommand{\eos}{\end{obsers}}
\newcommand{\bpf}{\begin{pf}}
\newcommand{\epf}{\end{pf}}
\newcommand{\ba}{\begin{array}}
\newcommand{\ea}{\end{array}}
\newcommand{\beq}{\begin{eqnarray}}
\newcommand{\beqq}{\begin{eqnarray*}}
\newcommand{\eeq}{\end{eqnarray}}
\newcommand{\eeqq}{\end{eqnarray*}}

\begin{document}
\bibliographystyle{amsplain}

\title{On Lipschitz continuity of solutions of hyperbolic Poisson's equation}

\author{Jiaolong Chen}
\address{Jiaolong Chen, Department of Mathematics,
Shantou University, Shantou, Guangdong 515063, People's Republic of China}
\email{jiaolongchen@sina.com}

\author{Manzi Huang}
\address{Manzi Huang, Department of Mathematics,
Hunan Normal University, Changsha, Hunan 410081, People's Republic of China}
\email{mzhuang79@163.com}

\author{Antti Rasila}
\address{Antti Rasila,
Department of Mathematics and Systems Analysis, Aalto University, P. O. Box 11100, FI-00076 Aalto, Finland.}
\email{antti.rasila@iki.fi}

\author{Xiantao Wang${}^{~\mathbf{*}}$}
\address{Xiantao Wang, Department of Mathematics,
Shantou University, Shantou, Guangdong 515063, People's Republic of China}
\email{xtwang@stu.edu.cn}

\keywords{Lipschitz continuity, hyperbolic Poisson's equation, Poisson kernel, Green function.\\
$^{\mathbf{*}}$Corresponding author}

%\dedicatory{}
\begin{abstract} In this paper, we investigate solutions of the hyperbolic Poisson equation $\Delta_{h}u(x)=\psi(x)$, where  $\psi\in L^{\infty}(\mathbb{B}^{n}, \IR^n)$ and
\[
\Delta_{h}u(x)= (1-|x|^2)^2\Delta u(x)+2(n-2)(1-|x|^2)\sum_{i=1}^{n} x_{i} \frac{\partial u}{\partial x_{i}}(x)
\]
is the hyperbolic Laplace operator in the $n$-dimensional space $\mathbb{R}^n$ for $n\ge 2$. We show that if $n\geq 3$ and $u\in C^{2}(\mathbb{B}^{n},\IR^n) \cap C(\overline{\mathbb{B}^{n}},\IR^n )$ is a solution to the hyperbolic Poisson equation, then it has the representation $u=P_{h}[\phi]-G_{ h}[\psi]$ provided that $u\mid_{\mathbb{S}^{n-1}}=\phi$ and $\int_{\mathbb{B}^{n}}(1-|x|^{2})^{n-1} |\psi(x)|\,d\tau(x)<\infty$. Here $P_{h}$ and $G_{h}$ denote Poisson and Green integrals with respect to $\Delta_{h}$, respectively. Furthermore, we prove that functions of the form $u=P_{h}[\phi]-G_{h}[\psi]$ are Lipschitz continuous.
\end{abstract}

\thanks{The research was partly supported by NSF of China (No. 11571216)}

\maketitle \pagestyle{myheadings} \markboth{Jiaolong Chen, Manzi Huang, Antti Rasila and Xiantao Wang}{On Lipschitz continuity of solutions of hyperbolic Poisson's equation}

\section{Introduction and main results}\label{csw-sec0}

For $n\ge 2$, let $\mathbb{B}^{n}(x_{0}, r)=\{x\in\mathbb{R}^{n}:|x-x_{0}|<r\}$, $\overline{\mathbb{B}^{n}}(x_{0}, r)=\{x\in\mathbb{R}^{n}:|x-x_{0}|\leq r\}$ and  $\mathbb{S}^{n-1}(x_{0}, r)=\partial\mathbb{B}^{n}(x_{0}, r)$. We write $\mathbb{B}^{n}=\mathbb{B}^{n}(0, 1)$ and $\mathbb{S}^{n-1}=\mathbb{S}^{n-1}(0, 1)$.

Let $L_{1}$, $L_{2}$ be two constants and $\Omega\subset\IR^n$ a domain. Then a mapping $f\colon  \Omega\to \mathbb{R}^{n}$ is said to be {\it  $L_{1}$-Lipschitz}  if $|f(x)-f(y)|\leq L_{1}|x-y|$ for all $x,y\in\Omega$, and {\it $L_{2}$-co-Lipschitz} if $|f(x)-f(y)|\geq L_{2}|x-y|$ for all $x,y\in\Omega$. If $f$ is both $L_{1}$-Lipschitz and $L_{2}$-co-Lipschitz for constants $L_{1}$ and $L_{2}$, then $f$ is called {\it bi-Lipschitz.}

In \cite{kapa2011}, Kalaj and Pavlovi\'{c} studied the bi-Lipschitz continuity of quasiconformal self-mappings of the unit disk $\ID=\mathbb{B}^{2}$ satisfying the Poisson's equation $\Delta u=\psi$, where $\Delta$ is the
usual Laplacian in $\mathbb{R}^n$.
See \cite{clw, chen, K1, KM0, KM, KP0, LCW, P} and references therein for further discussions along this line in the plane.

In \cite{AKM2008}, Arsenovi\'c et al. showed that the Lipschitz continuity of $\phi\colon \mathbb{S}^{n-1}\rightarrow\mathbb{R}^{n}$ implies the Lipschitz continuity of its harmonic extension
$P[\phi]\colon \mathbb{B}^{n}\rightarrow\mathbb{R}^{n}$ provided that $P[\phi]$ is a $K$-quasiregular mapping.
Here $P$ is the usual Poisson kernel with respect to $\Delta$, i.e.
\[
P[\phi](\eta) = \int_{\IS^{n-1}}P(\eta,\xi)\phi(\xi)\,d\sigma(\xi)\;\;\mbox{and}\;\; P(\eta,\xi)=\frac{1-|\eta|^2}{|\eta-\xi|^n},
\]
where $\eta\in \mathbb{B}^n$ and $\sigma$ is the $(n-1)$-dimensional Lebesgue measure normalized so that $\sigma(\IS^{n-1})=1$.
Moreover, Kalaj \cite{kalaj2008} also proved the Lipschitz continuity of $P[\phi]\colon \mathbb{B}^{n}\rightarrow\mathbb{B}^{n}$ under an additional assumption that it is a $K$-quasiconformal harmonic mapping with $P[\phi](0)=0$ and $\phi\in C^{1,\alpha}$ for some $\alpha \in (0,1]$. Later, in \cite{Kalaj2011}, Kalaj proved that $K$-quasiconformal mappings of $\mathbb{B}^{n}$ onto itself are Lipschitz, provided that they satisfy the Poisson equation  $\Delta u=\psi$ with $\psi\in L^{\infty}(\overline{\mathbb{B}^{n}}, \IR^n)$ and $u(0)=0$.

\subsection{Main results}
The purpose of this paper is to consider results of the above type for solutions of the {\it hyperbolic Laplace equation}.

\bdefe
 A function $u\in C^{2}(\mathbb{B}^{n}, \IR^n)$ $(n\ge 2)$ is said to be {\it hyperbolic harmonic} \cite{reka, sto1999, sto2016, sto2012} if it satisfies the hyperbolic Laplace equation
\be\label{eq00.002.1000}
\Delta_{h}u (x)=(1-|x|^2)^2\Delta u(x)+2(n-2)(1-|x|^2)\sum_{i=1}^{n} x_{i} \frac{\partial u}{\partial x_{i}}(x)=0.
\ee
 \edefe
Obviously, for $n=2$, hyperbolic harmonic functions coincide with harmonic functions. See \cite{clsh,du} for the properties of harmonic mappings. Also, see $\S$\ref{HyperP} below for more properties of $\Delta_{h}$.

It is well known that if $u$ satisfies the conditions: $(1)$ $\Delta u=\psi$ which is continuous in $\mathbb{B}^{n}$ with $n\ge 2$, and $(2)$ $u\mid_{\mathbb{S}^{n-1}}=\phi$ which is bounded and integrable in $\mathbb{S}^{n-1}$,
then (cf. \cite[p. 118-119]{LH} or \cite{Kalaj2011, kapa2011, kavu2015})
\[
u=P[\phi]-G[\psi]\;\;\mbox{and}\;\; G[\psi](\eta) = \int_{\IB^n}G(\eta,\xi)\psi(\xi)\,dV(\xi),
\]
where $V$ is the $n$-dimensional Lebesgue volume measure and $G(\eta,\xi)$, $\eta\neq \xi$, is the usual Green function \cite{Kalaj2011, kapa2011, kavu2015}, i.e.
\[
G(\eta,\xi)=
\left\{
\begin{array}{rl}
\frac{1}{2\pi}\log \left|\frac{1-\eta\overline{\xi}}{\eta-\xi}\right|, & \text{ for }n=2,\\
\frac{1}{(n-2)\omega_{n-1}}\Big(|\eta-\xi|^{2-n}-\big|\xi|\eta|-\eta/|\eta|\big|^{2-n}\Big), & \text{ for }n\geq3.
\end{array}\right.
\]
Here $\omega_{n-1} = 2\pi^{n/2}/\Gamma(n/2)$ is the $(n-1)$-dimensional surface area of $\IS^{n-1}$ and $\Gamma$ is the Gamma function (see e.g.  \cite[p. 61]{ahlf} or \cite[Appendix A]{axbo}).

The first aim of this paper is to establish the counterpart of the above result to the solutions to the Dirichlet problem:
 \be\label{eq-1.1}
 \left\{\begin{array}{rl}
\Delta_{h}u(x)=\psi(x),&x\in\mathbb{B}^{n}, \\
u(\xi)=\phi(\xi),&\xi\in\mathbb{S}^{n-1},
\end{array}\right.
\ee
where $\psi\in L^{\infty}(\mathbb{B}^{n}, \IR^n)$ and $\phi\in L^{\infty}(\mathbb{S}^{n-1}, \IR^n)$.

Denote by $\tau$ the {\it M\"{o}bius invariant measure} in $\mathbb{B}^{n}$, which is given by
 $$d\tau(x)=\frac{d\nu(x)}{(1-|x|^2)^n},$$
where $\nu$ is the $n$-dimensional Lebesgue volume measure normalized so that $\nu(\IB^n)=1$.
Our result is as follows:

\begin{thm}\label{thm-3}
Suppose $u\in C^{2}(\mathbb{B}^{n},\IR^n) \cap C(\overline{\mathbb{B}^{n}},\IR^n)$, $n\geq3$ and
\[
\int_{\mathbb{B}^{n}}(1-|x|^{2})^{n-1}|\psi(x)|\,d\tau(x)\leq \mu_{1},
\]
where $\mu_{1}\geq 0$ is a constant. If $u$ satisfies \eqref{eq-1.1}, then
\be\label{jl-1}u=P_{h}[\phi]-G_{ h}[\psi].\ee

\end{thm}

%\noindent
Here $P_{h}[\phi]$  and $G_{ h}[\psi]$  denote the Poisson integral of $\phi$ and the Green integral of $\psi$, with respect to $\Delta_h$, respectively (See \eqref{eq20.20.10} and \eqref{eq20.20.20} below for the details).

The second aim of this paper is to establish the Lipschitz continuity of the mappings $u$ of the form \eqref{jl-1}. More precisely, we have the following.

\begin{thm}\label{thm-4}
Let $n\ge 3$. Suppose
 \ben
\item
$u\in C^{2}(\mathbb{B}^{n}, \IR^n) \cap C(\overline{\mathbb{B}^{n}},\IR^n )$ is of the form \eqref{jl-1};
 \item
there is a constant $L\geq 0$ such that $|\phi(\xi)-\phi(\eta)|\leq L|\xi-\eta|$ for all $\xi,\eta\in \mathbb{S}^{n-1}$;
 \item
there is a constant $M\geq 0$ such that $|\psi(x)|\leq M(1-|x|^{2}) $ for all $x\in \mathbb{B}^{n}$.
\een
Then there exists a constant $N=N(n, L,M)$ such that for $x$, $y\in \mathbb{B}^{n}$,
  $$|u(x)-u(y)| \leq  N |x-y|.$$
 \end{thm}

\br
In Section \ref{example-sec}, we give an example to show that the assumption ``$n\geq 3$" in Theorem \ref{thm-4} is necessary.
\er

In fact, Theorem \ref{thm-4} follows from more general, albeit technical, results on Lipschitz continuity of $P_{h}[\phi]$ and $G_{h}[\psi]$, which we shall discuss next.

\subsection{$\omega$-Lipschitz continuity}

 A continuous increasing function $\omega\colon [0,\infty)\rightarrow[0,\infty)$ with $\omega(0)=0$ is called a {\it majorant} if $\omega(t)/t$ is non-increasing for $t>0$.
 Given a subset $\Omega$ of $\mathbb{R}^{n}$, a function $f:\Omega\rightarrow\mathbb{R}^{n}$ is said to be {\it $\omega$-Lipschitz continuous} or belong to the {\it Lipschitz space $\Lambda_{\omega}(\Omega)$} if there is a positive constant $C$ such that
\be\label{eq21.21.21}|f(x)-f(y)|\leq C \omega(|x-y|)\ee
for all $x,y\in \Omega$ (cf. \cite{chenshao2015, dyak1997, dyak2004, Pa2004, Pa2007}).
For some $\rho_{0}>0$ and $0<\rho<\rho_{0}$, a majorant $\omega$ is called {\it fast} if
% \be\label{eq31.31.31}
$$ \int_{0}^{\rho}\frac{\omega(t)}{t}dt\leq C\omega(\rho).$$

Let $\Omega$ be a proper subdomain of $\mathbb{R}^{n}$. We say that a function $f\colon \Omega\rightarrow \mathbb{R}^{n}$ belongs to the {\it local Lipschitz space ${\rm loc}\Lambda_{w}(\Omega)$} if \eqref{eq21.21.21} holds, whenever $x\in \Omega$ and $|x-y|<\frac{1}{2}\delta_{\Omega}(x)$,
where $C$ is a positive constant and $\delta_{\Omega}(x)$ denotes the Euclidean distance from $x$ to the boundary $\partial \Omega$ of $\Omega$.

A domain $\Omega\subset \IR^n$ is said to be a {\it $\Lambda_{w}$-extension domain} if $\Lambda_{w}(\Omega)={\rm loc}\Lambda_{w}(\Omega)$.
In \cite{Lappalainen}, Lappalainen proved that $\Omega$ is a $\Lambda_{w}$-extension domain if and only if each pair of points $x,y\in \Omega$ can be joined by a rectifiable curve $\gamma \subset \Omega$ satisfying
 \be\label{eq21.21.221} \int_{\gamma}\frac{ \omega(\delta_{\Omega}(\eta))}{\delta_{\Omega}(\eta)}ds(\eta)\leq C \omega(|x-y|)\ee
with some fixed positive constant $C=C(\Omega,\omega)$ which means that the constant $C$ depends only on the quantities $\Omega$ and $\omega$, where $ds$ is the length measure on $\gamma$.
Furthermore, we know from \cite[Theorem 4.12]{Lappalainen} that $\Lambda_{w}$-extension domains exist for fast majorants $\omega$ only. Conversely, if
$\omega$ is fast, then the class of $\Lambda_{\omega}$-extension domains is fairly large and contains all bounded uniform domains.

\br
Recall that a domain $\Omega$ is said to be {\it uniform} if there is a constant $C$ such that
each pair of points $x_{1}$ and $x_{2}$ in $\Omega$ can be joined by a rectifiable curve $\gamma\subset\Omega$ satisfying
\begin{center}
$\ell(\gamma)\leq C|x_{1}-x_{2}|$ and
$\min \{\ell(\gamma[x_{1}, x]),\; \ell(\gamma[x_{2}, x])\}\leq C \delta_{\Omega}(x)$
\end{center}
 for all $x\in \gamma$.
Here $\ell(\gamma)$ denotes the length of $\gamma$ and $\gamma[x_{i}, x])$ is the subarc of $\gamma$ with endpoints $x_i$ and $x$, where $i=1$, $2$.
It is known that $\mathbb{B}^{n}$ is a uniform domain, and hence a $\Lambda_{w}$-extension domain for a fast $\omega$ \cite[Section 1]{dyak2004}.
\er

The next two results establish $\omega$-Lipschitz continuity of $P_{h}[\phi]$ and $G_{h}[\psi]$:

\begin{thm}\label{thm-11}
Suppose $n\geq 3$, $\phi\colon \mathbb{S}^{n-1}\rightarrow \mathbb{R}^{n}$ and $|\phi(\xi)-\phi(\eta)|\leq \omega(|\xi-\eta|)$ for all $\xi,\eta\in \mathbb{S}^{n-1}$, where $\omega$ is a fast majorant.
 Then, for $x,y\in \mathbb{B}^{n}$, $$|\Phi(x)-\Phi(y)|\leq C \alpha_{0} \omega(|x-y|),$$
 where $\Phi=P_{h}[\phi]$ and $\alpha_{0}=\alpha_{0}(n)$ and $C=C(\mathbb{B}^{n},\omega)$
 is the same constant as in \eqref{eq21.21.221}.
\end{thm}

\begin{thm}\label{thm-2}
Suppose $n\geq 3$, $\psi\in C (\mathbb{B}^{n},\IR^n)$ and $|\psi(x)|\leq M(1-|x|^{2}) $ for $x\in \mathbb{B}^{n}$, where $M$ is a constant.
 Then, for $x,y\in \mathbb{B}^{n}$,
$$|\Psi(x)-\Psi(y)|\leq \beta_{0}|x-y|,$$
 where $\Psi=G_{ h}[\psi]$ and $\beta_{0}=\beta_{0}(n,M)$ is a constant.
\end{thm}

This paper is organized as follows. In Section \ref{csw-sec1}, some necessary terminology and known results are introduced, and several preliminary results are proved.
In Section \ref{sec-2}, we present the proof of Theorem \ref{thm-3}. In Section \ref{LipP}, we show Theorem \ref{thm-11}.
In Section \ref{LipG}, we prove Theorem \ref{thm-2} and Theorem \ref{thm-4}.
Finally, in Section \ref{example-sec}, we construct two examples to illustrate
the necessity of the requirement $n\geq 3$ in Theorem \ref{thm-4},
and the existence of hyperbolic harmonic mappings which satisfy the conditions in Theorem \ref{thm-11} but not $K$-quasiregular.

\section{Preliminaries}\label{csw-sec1}

In this section, we recall some necessary terminology and results.

\subsection{Matrix notations}

For natural number $n$, let
$$A=\big(a_{ij}\big)_{n\times n}\in \mathbb{R}^{n\times n}. $$

For $A\in \mathbb{R}^{n\times n},$
denote by $\|A\|$ the matrix norm $\|A\|=\sup\{|Ax|:\; x\in \mathbb{R}^{n}, |x|=1\}$, and $l(A)$ the matrix function $l(A)=\inf\{|Ax|:\; x\in \mathbb{R}^{n}, |x|=1\}.$

For a domain $\Omega\subset\mathbb{R}^{n}$, let $u=(u_{1},\ldots,u_{n})\colon \Omega\rightarrow \IR^n$ be a function that has all partial derivatives at $x=(x_{1},\ldots,x_{n})$ on $\Omega$. Then $Du$ denotes the usual Jacobian matrix
\be\label{ja}
Du=\left(
    \begin{array}{cccc}
      \frac{\partial u_{1}}{\partial x_{1}}    & \cdots & \frac{\partial u_{1}}{\partial x_{n}} \\
     \vdots & \ddots  & \vdots \\
    \frac{\partial u_{n}}{\partial x_{1}}  & \cdots & \frac{\partial u_{n}}{\partial x_{n}}\\
    \end{array}
  \right) = \big(\nabla u_1 \cdots \nabla u_n\big)^T,
\ee
where $T$ is the transpose and the gradients $\nabla u_j$ are understood as column vectors. If $Du$ is a nonsingular matrix, then the eigenvalues $\lambda_j^2$ of the (symmetric and positive definite) matrix $Du\times Du^T$ are real, and they can be ordered so that $0<\lambda_{1}^{2}\leq \lambda_{2}^{2}\leq\ldots\leq\lambda_{n}^{2} $. Then
$|J_{u}|=\prod_{k=1}^{n}\lambda_{k}$, $l(D u)=\lambda_{1}$ and $\|D u\|=\lambda_{n}$,
where $J_{u}$ denotes the Jacobian of $u$.

\subsection{Spherical coordinate transformation}
Let $Q=(\xi_{1},\ldots,\xi_{n}):K^{n-1}\rightarrow \mathbb{S}^{n-1}$ be the following spherical coordinate transformation \cite{kalaj2008}:
\be\label{Ttrans}
\begin{array}{rcl}
\xi_{1} &=&\cos\theta_{1},\\
\xi_{2}&=&\sin\theta_{1}\cos\theta_{2},\\
&\vdots &\\
\xi_{n-1} &=& \sin\theta_{1} \sin\theta_{2} \ldots  \sin\theta_{n-2}\cos\theta_{n-1},\\
\xi_{n} &=& \sin\theta_{1} \sin\theta_{2} \ldots \sin\theta_{n-2}\sin\theta_{n-1}
.
\end{array}
\ee
Here $K^{n-1}= [0,\underbrace{\pi]\times\ldots\times [0}_{n-2},\pi] \times[0,2\pi]$. Note that
\be\label{eq1.0.0} J_{Q}(\theta_{1},\ldots,\theta_{n-1})=\underbrace{\sin^{n-2}\theta_{1}\ldots \sin \theta_{n-2}}_{n-2}.\ee

For an integrable function $f$ in $\mathbb{B}^{n}$, by letting $x=\rho\xi$ with $\rho=|x|$, we have
\be\label{eq0.0cj.0}
 \int_{\mathbb{B}^{n}(0, r)}f(x)d\nu(x)=n\int_{0}^{r}\rho^{n-1}d\rho\int_{\mathbb{S}^{n-1}}f(\rho\xi) d\sigma(\xi),
\ee
where $$d\sigma(\xi)=\frac{1}{\omega_{n-1}} J_{Q}(\theta_{1},\ldots,\theta_{n-1})  d\theta_{1}\ldots d\theta_{n-1}$$
(see, e.g. \cite{Kalaj2011,sto1999,zhuke}).

\subsection{$K$-quasiregular}
A sense-preserving continuous mapping $u:$ $\Omega\to G$ between two open subsets $\Omega$ and $G$ of the Euclidean space $\mathbb{R}^{n}$ will be called a {\it $K$-quasiregular mapping} ($K\geq 1$)
if\begin{enumerate}
\item
$u$ is an absolutely continuous mapping in almost every segment parallel to some of the coordinate axes and there exist the partial derivatives which are locally $L^{n}$ integrable mappings on $\Omega$. We will write $u\in ACL^{n}$ , and
\item
$u$ satisfies the condition $||Du(x)||^{n}/K \leq J_{u}(x)\leq Kl\big(D u(x)\big)^{n}$ at $x$ almost everywhere in $\Omega$.
\end{enumerate}

In particular, $u$ is called $K$-quasiconformal if $u$ is a $K$-quasiregular homeomorphism (cf. \cite{kalaj2008, rick1993}).

\subsection{Hypergeometric functions}
Let $F$ be the hypergeometric function
\be\label{eq2.2.2}F(a,b;c;s)=\sum_{k=0}^{\infty}\frac{(a)_{k}(b)_{k}}{k!(c)_{k}}s^{k},\ee
where $a,b,c\in \mathbb{R}$, $c$ is neither zero nor a negative integer,
$(a)_{k}$ denotes the Pochhammer symbol with $(a)_{0}=1$
and $(a)_{k}=a(a+1)\ldots (a+k-1)$ $(k\in \mathbb{N})$.
If $a$ is not a negative integer, then
$$(a)_{k}=\Gamma(a+k)/\Gamma(a).$$
If $c-a-b>0$, then the series \eqref{eq2.2.2} converges absolutely for all $|s|\leq 1$ (cf. \cite[Section 31]{rain}).

Let $t>1$, $  k\in \mathbb{R}$ and $r\in (-1,1)$. Ren and K\"{a}hler \cite[Lemma 2.2]{reka} proved that
\be\label{eq2.2.4}\int_{-1}^{1}\frac{(1-s^2)^{(t-3)/2}}{(1-2rs+r^2)^{k}}ds=\frac{\Gamma(\frac{t-1}{2})
\Gamma(\frac{1}{2})}{\Gamma(\frac{t}{2})}F\left(k,k+1-\frac{t}{2};\frac{t}{2};r^{2}\right).\ee

The following lemmas will be useful in the proof of Theorem \ref{thm-4}.
\begin{lem}\label{lem-4.5} Let
 $$f_{n}(s)=\sum_{k=0}^{\infty}\frac{(a)_{k} (\frac{b-n}{2})_{k}}{(\frac{n}{2})_{k}k!}\frac{s^k}{k+c}.$$
Suppose $b$, $n\in \mathbb{N}^{+},$ $a>0$, $c>0$ and $n-a-\frac{b}{2}>0$. Then there is a constant $\mu_2\geq 0$ such that for all $s\in [0,1]$ and all $n\geq b$,
$$|f_{n}(s)|\leq \mu_2,$$
where $\mu_2=\mu_2(n,a,b,c)$.
\end{lem}
\bpf Obviously, we only need to consider the case where $b$ is even since the proof of the case $b$ being odd is similar.
To finish the proof, we consider the following two possibilities.

 \begin{case}\label{jl131-3-even} $n$ is even.\end{case}
Under this assumption, we easily see from  $n\geq b$ that
$$\Big(\frac{b-n}{2}\Big)_{k}=\Big(\frac{b-n}{2}\Big)\Big(\frac{b-n}{2}+1\Big)\ldots\Big(\frac{b-n}{2}+k-1\Big)=0$$ for all $k\geq \frac{n-b+2}{2}$,
and hence $f_{n} $ is a polynomial,
where $$f_{n}(s)= \sum_{k=0}^{\frac{n-b+2}{2}} \frac{(a)_{k} (\frac{b-n}{2})_{k}}{(\frac{n}{2})_{k}k!}\cdot\frac{s^k}{k+c}.$$
Hence, for all $s\in [0,1]$,
\be\label{jl-4}|f_{n}(s)|\leq \mu',
\text{ where }\mu'=\sum_{k=0}^{\frac{n-b+2}{2}} \frac{(a)_{k} |(\frac{b-n}{2})_{k}|}{(\frac{n}{2})_{k}k!}\frac{1}{k+c}.
\ee

\begin{case}\label{jl131-3-odd} $n$ is odd.\end{case}
In this case, we separate $f_{n} $ into two parts: $f_{n} =f_{n_{1}} +f_{n_{2}}$,
%\be\label{eq04.1.0}
where $$f_{n_{1}}(s)=\sum_{k=0}^{ \frac{n-b+1}{2}}\frac{(a)_{k} (\frac{b-n}{2})_{k}}{(\frac{n}{2})_{k}k!}\frac{s^k}{k+c} $$
 and
 $$f_{n_{2}}(s)=\sum_{k= \frac{n-b+3}{2} }^{\infty}\frac{(a)_{k} (\frac{b-n}{2})_{k}}{(\frac{n}{2})_{k}k!}\frac{s^k}{k+c}.$$

Since $f_{n_{1}}$ is continuous in $[0,1]$, obviously,  for $s\in [0,1]$,
\be\label{jl-5} |f_{n_{1}}(s)|\leq \mu'',
\text{ where }\mu''=\sum_{k=0}^{ \frac{n-b+1}{2}}\frac{(a)_{k} |(\frac{b-n}{2})_{k}|}{(\frac{n}{2})_{k}k!}\frac{1}{k+c}.
\ee

Next, we estimate $f_{n_{2}}$.
Since $\frac{b-n}{2}+t-1>0$ for $t\geq \frac{n-b+3}{2}$ and $\frac{b-n}{2}+t-1<0$ for $t\leq \frac{n-b+1}{2}$, we obtain that
$$f_{n_{2}}(s)=(-1)^{\frac{n-b+1}{2}}\sum_{k= \frac{n-b+3}{2} }^{\infty}\frac{(a)_{k} |(\frac{b-n}{2})_{k}|}{(\frac{n}{2})_{k}k!}\frac{s^k}{k+c}.$$
We leave the estimate on $f_{n_{2}}$ for a moment and prove the following claim.
\bcl\label{tue-1}
 Let
$$g_{a, b, n}(s)=\sum_{ k= \frac{n-b+3}{2} }^{\infty}\frac{(a)_{k} (\frac{b-n}{2})_{k}}{(\frac{n}{2})_{k}k!} s^k.$$ Then $g_{a, b, n}$ is continuous in $[0,1]$.
\ecl
To prove the continuity of $g_{a, b, n}$ in $[0, 1]$, it suffices to check the uniform convergence of  $g_{a, b, n}$ in $[0, 1]$. Since
$$g_{a, b, n}(s)=(-1)^{\frac{n-b+1}{2}}\sum_{k= \frac{n-b+3}{2} }^{\infty}\frac{(a)_{k} |(\frac{b-n}{2})_{k}|}{(\frac{n}{2})_{k}k!} s^k,$$
obviously, we only need to demonstrate the boundedness of
$$\sum_{k= \frac{n-b+3}{2} }^{\infty}\frac{(a)_{k} |(\frac{b-n}{2})_{k}|}{(\frac{n}{2})_{k}k!}.$$

This easily follows from the following two facts:
\ben
\item  It follows from the assumption $``\frac{n}{2}-a-\frac{b-n}{2}>0"$ and \cite[Section 31]{rain} that
$$F\left(a,\frac{b-n}{2};\frac{n}{2};s\right)=\sum_{ k= 0}^{\infty}\frac{(a)_{k} (\frac{b-n}{2})_{k}}{(\frac{n}{2})_{k}k!} s^k$$ is bounded in $[0,1]$.
\item
$\sum_{ k= 0}^{k= \frac{n-b+1}{2}}\frac{(a)_{k} (\frac{b-n}{2})_{k}}{(\frac{n}{2})_{k}k!} s^k$ is continuous in $[0, 1].$
\een\medskip

Now, we continue the estimate on $f_{n_2}$.
Let
$$\| g_{a, b, n} \|_{\infty}=\max\{|g_{a, b, n}(s)|: \; s\in [0, 1]\}.$$
Then Claim \ref{tue-1} guarantees that $\| g_{a, b, n} \|_{\infty}$ is finite.
It follows that for all $s\in [0,1]$,
\be\label{jl-6} | f_{n_{2}}(s)|\leq |g_{a, b, n}(s)|\leq \| g_{a, b, n} \|_{\infty}.\ee
By taking $$\mu_2=\max\{\mu', \mu''+\| g_{a, b, n} \|_{\infty}\},$$
the lemma follows from \eqref{jl-4}, \eqref{jl-5} and \eqref{jl-6}.    \epf

\begin{lem}\label{jl-30}
Let
$$I_{0}(s)=\int_{0}^{1}t^{m} F\left(a,\frac{b-n}{2};\frac{n}{2};ts \right)dt.$$
Suppose $b$, $n \in \mathbb{N}^{+}$, $a>0$ and $m> -1$. Then for all $s\in[0,1)$,
$$I_{0}(s)=\sum_{k=0}^{\infty}\frac{(a)_{k} (\frac{b-n}{2})_{k}}{(\frac{n}{2})_{k}k!}\frac{s^k}{k+m+1}.$$
\end{lem}
\bpf  Obviously,
$$ \int_{0}^{1}t^{m} F\left(a,\frac{b-n}{2};\frac{n}{2};ts \right)dt=
\int_{0}^{1} \sum_{k=0}^{\infty}\frac{(a)_{k} (\frac{b-n}{2})_{k}}{(\frac{n}{2})_{k}k!}s^k t^{k+m} dt,$$
and the convergence radius of the series $\sum_{k=0}^{\infty}\frac{(a)_{k} (\frac{b-n}{2})_{k}}{(\frac{n}{2})_{k}k!}s^k t^{k+m} $
is $1/s$ for $s\in [0,1)$. Hence we have
%\begin{eqnarray*}
$$I_{0}(s)
= \sum_{k=0}^{\infty}\frac{(a)_{k} (\frac{b-n}{2})_{k}}{(\frac{n}{2})_{k}k!}s^k\int_{0}^{1}t^{k+m}  dt
% \\ \nonumber
=\sum_{k=0}^{\infty}\frac{(a)_{k} (\frac{b-n}{2})_{k}}{(\frac{n}{2})_{k}k!}\frac{s^k}{k+m+1},$$
 %\end{eqnarray*}
as required.
\epf

\begin{lem}\label{lem-4.2}
Suppose $b$, $n \in \mathbb{N}^{+}$, $a>0,$ $n-a-\frac{b}{2}>0$, $n\geq b$ and $m> -1$. Then for all $s\in[0,1)$,
$$|I_{0}|\leq \mu_{2,1},$$
where $\mu_{2,1}=\mu_{2}(n,a,b,m+1)$ and $I_{0}$ are defined in Lemmas \ref{lem-4.5} and \ref{jl-30}, respectively.
\end{lem}
\bpf This lemma easily follows from Lemma \ref{lem-4.5} and Lemma \ref{jl-30}.
\epf

\subsection{M\"{o}bius transformations}

 For any $x,$ $y\in \mathbb{R}^{n}$, we denote the inner product $\sum_{k=1}^{n}x_{k}a_{k}$ by $\langle x,a\rangle$.
 Let $x=|x|x'$ and $y=|y|y'$. Then the symmetry lemma (see e.g. \cite{ahlf} or \cite{axbo, reka}) shows that
 $$\big||y|x-y'\big|=\big||x|y-x'\big|.$$
In the following, we denote $[x,y]=\big||x|y-x'\big|$. Obviously, $[x,y]=[y,x]$.

 For any $a\in \mathbb{B}^{n}$, let
%\be\label{eq1.0.1}
$$\varphi_{a}(x)=\frac{|x-a|^{2}a-(1-|a|^{2})(x-a)}{[a,x]^{2}}$$
 in $\mathbb{B}^{n}$.
 Then $\varphi_{a}$ is a M\"{o}bius transformation of $\mathbb{R}^{n}$ mapping $\overline{\mathbb{B}^{n}}$ onto $\overline{\mathbb{B}^{n}}$
 with $\varphi_{a}(a)=0$, $ \varphi_{a}(0)=a$ and $\varphi_{a}(\varphi_{a}(x))=x$  \cite{sto2012}.
It follows from Equations (2.4) and (2.6), Theorem 3.4(a) and Chapter 5 in \cite{sto1999}, together with \cite[Equation (2.4)]{reka}, that
\beq \label{eq1.0.2}[a,x]^{2}&=&|x-a|^{2}+(1-|x|^2)(1-|a|^2)\\  \nonumber
&=&1+|a|^{2}|x|^{2}-2|a| |x | \left\langle\frac{x}{|x|},\frac{a}{|a|}\right\rangle, \eeq
\be\label{wed-1}
 1-|x|  \leq [a,x] < 2,
\ee
\be\label{eq1.0.3} |\varphi_{a}(x)|=|\varphi_{x}(a)|=\frac{|x-a|}{[a,x]},\;\;\; 1-|\varphi_{a}(x)|^2=\frac{(1-|x|^2)(1-|a|^2)}{[a,x]^2},\ee
\be\label{eq1.000.3}J_{\varphi_{a}}(x)=\frac{(1-|a|^2)^{n}}{[a,x]^{2n}}
\ee
and
\beq \label{eq1.0.5}\;\;\;\;\;\;\frac{\partial}{\partial x_{k}}|\varphi_{a}(x)| &=&\frac{\partial}{\partial x_{k}}|\varphi_{x}(a)|\\  \nonumber
&=&\frac{[a,x]^2(x_{k}-a_{k})-|a-x|^{2}(x_{k}-a_{k})+|a-x|^{2}(1-|a|^2)x_{k}}{|a-x| \cdot[a,x]^3}.
 \eeq

Elementary calculations lead to
\be\label{eq1.0.4}  |a-\varphi_{a}(x)| =\frac{(1-|a|^2) |x| }{[a,x] }\;\;\mbox{and}\;\;\; [a,\varphi_{a}(x)] =\frac{1-|a|^2 }{[a,x] }.\ee

 We denote by $\mathcal{M}(\mathbb{B}^{n})$ the set of all M\"{o}bius transformations in $\mathbb{B}^{n}$. By \cite[Theorem 2.1]{sto1999} or
 \cite[Section 2]{sto2012},
 if $\varphi\in \mathcal{M}(\mathbb{B}^{n})$, then there exist $a\in \mathbb{B}^{n}$ and
 an orthogonal transformation $A$ such that $$\varphi(x)=A\varphi_{a}(x).$$
For more information about the M\"{o}bius transformations in $\mathbb{B}^{n}$, see e.g. \cite[Chapter 2]{ahlf}, \cite{beardon} or \cite[Chapter 1]{matti}.

\subsection{Hyperbolic Poisson's equation }\label{HyperP}
In terms of the mapping $ \varphi_{a}$, the {\it hyperbolic metric} $d_{h}$ in $\mathbb{B}^{n}$ is given by
$$d_{h}(a,b)=\log \left(\frac{1+|\varphi_{a}(b)|}{1-|\varphi_{a}(b)|}\right)$$
for all $ a,b\in \mathbb{B}^{n}$.

For all $\varphi \in \mathcal{M}(\mathbb{B}^{n})$, by the definition of $\Delta_{h}$, we have the following M\"obius invariance property \cite[Section 2]{sto2012}:
  \be\label{eq0.0.1}\Delta_{h}( u\circ \varphi )=\Delta_{h} u\circ \varphi .\ee
Obviously,
\be\label{eq00.00.1000} \Delta_{h} u (x)= \Delta(u\circ \varphi_{x})(0).\ee
In fact, if \eqref{eq00.00.1000} holds for all $u\in C^2(\mathbb{B}^{n})$ and $x\in \mathbb{B}^{n}$, then we can show
that $\Delta_{h}$ has the representation \eqref{eq00.002.1000} (cf. \cite[Chapter 3]{sto1999} or \cite{sto2016}).

Let
 \be\label{eq0.0.01} g(r,t)=\frac{1}{n}\int_{r}^{t}\frac{(1-s^2)^{n-2}}{s^{n-1}}ds \;\;\;\mbox{and} \;\;\; g(r)=g(r,1),\ee
where $0\leq r<t<1$.
 It is well known that the Green's function $G_{h}(x,y)$ w.r.t. $\Delta_{h}$ is given by
 \be\label{eq0.0.210}G_{h}(x,y)=g(|\varphi_{x}(y)|)=\frac{1}{n}\int_{|\varphi_{x}(y)|}^{1}\frac{(1-s^2)^{n-2}}{s^{n-1}}ds\ee
 for all $x\neq y\in \mathbb{B}^{n}$.

We remark that in the complex plane $\mathbb{C}$, every M\"{o}bius transformation $\varphi$ mapping the unit disc $\mathbb{D}$ onto itself can be written as $\varphi(z)=e^{i\theta} \varphi_{w}(z)$, where $\varphi_{w}(z)=\frac{w-z}{1-\overline{w}z}$ for some $w$ in $\mathbb{D}$. Hence when $n=2$, by \eqref{eq0.0.210},
we get \cite{kapa2011}
\be\label{cgeq0.0.210} G_{h}(w,z)=g\big(|\varphi_{w}(z)|\big)= \frac{1}{2}\log\frac{|1-\overline{w}z|}{|w-z|}= \pi\cdot G(w,z),\ee
where $G$ is the usual Green function w.r.t. $\Delta$.\medskip

For function $g$ in \eqref{eq0.0.01}, we define
\be\label{sun-2}q(t)=\frac{t^{n-2}}{(1-t^{2})^{n-1}}g(t)
\ee
in $(0,1)$.
Since elementary calculations lead to $$\lim_{t\rightarrow 0^{+}}q(t)=\frac{1}{n(n-2)} \;\;\text{and}\;\;\lim_{t\rightarrow1^{-}}q(t)=\frac{1}{2n(n-1)},$$ we define $$q(0)=\frac{1}{n(n-2)}\;\; \mbox{and}\;\;q(1)=\frac{1}{2n(n-1)}.$$
Then we have
\blem\label{eq0.0.21} For $n\geq 3$,
$\frac{1}{2n(n-1)}
\leq q(t)\leq
\frac{1}{n(n-2)}$ in $[0, 1]$.
\elem
\bpf
We start with the following claim.
\bcl
For $n\geq 3$,
$\frac{1}{2n(n-1)}
< q(t)<
\frac{1}{n(n-2)}$ in $(0, 1)$.
\ecl
 For $t\in (0, 1]$, let
$$q_1(t)=g(t)-\frac{1}{n(n-2)}\cdot  \frac{(1-t^{2})^{n-1}}{t^{n-2}}$$
and
$$q_2(t)=g(t)-\frac{1}{2n(n-1)}\cdot  \frac{(1-t^{2})^{n-1}}{t^{n-2}}.$$
Then $q_1(t)$ is increasing and $q_2(t)$ is decreasing, respectively, in $(0,1)$. Since $q_1(1)=q_2(1)=0$, we see that
$$ \frac{1}{2n(n-1)}\cdot  \frac{(1-t^{2})^{n-1}}{t^{n-2}} \leq g(t)\leq \frac{1}{n(n-2)}\cdot  \frac{(1-t^{2})^{n-1}}{t^{n-2}}$$ in $(0, 1]$, which implies that the claim holds.\medskip

Now, the lemma easily follows from Claim \ref{tue-1} and \eqref{sun-2}.
 \epf

The {\it Poisson-Szeg\H{o} kernel} $P_{h}$ for $\Delta_{h}$ is given by
 \be\label{eqc.0.21}P_{h}(x,t)=\left(\frac{1-|x|^2}{|t-x|^2}\right)^{n-1},\ee
which satisfies (cf. \cite[Lemma 5.20]{sto1999} or \cite{sto2016})
$$ \int_{\mathbb{S}^{n-1}}P_{h}(x,t)d \sigma(t)=1,$$
 and for each $k\in \{1, 2, \ldots, n\}$,
\beq \label{eq1.1.0}
 &\;\;&\frac{\partial}{\partial x_{k}}P_{h}(x,t)=\frac{\partial}{\partial x_{k}} \left( \frac{1-|x|^{2}}{|t-x|^{2}}\right)^{n-1}\\ \nonumber
&=& -2(n-1)\frac{x_{k}|t-x|^{2}+(1-|x|^{2})(x_{k}-t_{k})}{|t-x|^{4}}\cdot\left( \frac{1-|x|^{2}}{|t-x|^{2}}\right)^{n-2},
\eeq
 where $(x,t)\in\mathbb{B}^{n}\times \mathbb{S}^{n-1}$.

If $\phi\in L^1(\mathbb{S}^{n-1}, \IR^n)$ $(n\geq 2)$, we define the Poisson-Szeg\H{o} integral or invariant Poisson integral of $\phi$
(cf. \cite{ahco, rogr, sto2016} or \cite[Definition 5.21]{sto1999}) by
\be\label{eq20.20.10}P_{h}[\phi](x)=\int_{\mathbb{S}^{n-1}}P_{h}(x,\xi)\phi(\xi)\,d\sigma(\xi).\ee
If $\psi$ satisfies the following conditions:
\begin{enumerate}
  \item For $n\geq 3$, $\psi\in C(\mathbb{B}^{n},\IR^n)$ and $\int_{\mathbb{B}^{n}}(1-|x|^{2})^{n-1} |\psi(x)|\,d\tau(x)<\infty$,
  \item For $n=2$, $\psi(z)=(1-|z|^2)^2\psi_{0}(z)$, where $\psi_{0}\in C(\mathbb{D},\mathbb{C})$,
\end{enumerate}
then we define a function $G_{ h}[\psi]$ by
 \beq\label{eq20.20.20}G_{ h}[\psi](x)&=&\int_{\mathbb{B}^{n}}G_{h}(x,y)\psi(y)\,d\tau (y)\\ \nonumber &=& \frac{1}{n}\int_{\mathbb{B}^{n}}\left[\psi(y)\int_{|\varphi_{x}(y)|}^{1}\frac{(1-s^2)^{n-2}}{s^{n-1}}ds\right]\,d\tau (y).\eeq
This function is called the invariant Green integral of $\psi$.

\br
 If $n=2$ and $\Delta u(z)\in C(\mathbb{D},\mathbb{C})$, then it follows from \eqref{cgeq0.0.210}, together with the facts
 $\Delta_{h}u(z)=(1-|z|^2)^2 \Delta u(z)$ and $d\tau(z)=\frac{1}{\pi} (1-|z|^{2})^{-2}dA(z)$,
 that
 \be\label{eqcj20.c20.c10}G_{h} [\Delta_{h} u](z) =\frac{1}{2\pi}\int_{ \mathbb{D}} \log \frac{|1-\overline{w}z|}{|w-z|}  \Delta_{h} u(z) \frac{dA(z)}{(1-|z|^2)^2}=G[\Delta u](z),\ee
 where $dA(re^{i \theta})=r\,dr\,d\theta$.

Furthermore, \eqref{eqc.0.21} implies that
 $P_{h} =P $ provided that $n=2$.
 Let $$\psi(z)=(1-|z|^{2})^{2}\psi_{0}(z),$$ where $\psi_{0}\in C(\mathbb{D},\mathbb{C})$.
 If $u$ satisfies $\Delta u=\psi_{0}$ in $\mathbb{D}$
 and $u|_{\mathbb{S}^{1}}=\phi\in L^{1}(\mathbb{S}^{1},\IC)$, then it follows from \cite{kapa2011}, \eqref{eq20.20.10} and \eqref{eqcj20.c20.c10} that
 \be\label{eqcg20.20.20}\nonumber  u=P[\phi]-G[\psi_{0}]=P_{h}[\phi]-G_{h}[\psi]. \ee
\er

We use $C^2_{c}(\mathbb{B}^{n})$ to denote the set of all twice continuous differentiable functions with compact support in $\mathbb{B}^{n}$.
Let us recall the following two results.

\begin{Thm}\label{Thm-D} {\rm(}\cite[Corollary 4.4]{sto1999} {\rm or} \cite{sto2016}{\rm)}
If $u\in C_{c}^{2}(\mathbb{B}^{n },\IR^n)$, then for all $x\in \mathbb{B}^{n }$,
$$u=-G_{ h}[\Delta_{h}u].$$
\end{Thm}

\begin{Thm}\label{Thm-B} {\rm(}\cite[Theorem 5.22]{sto1999} {\rm or} \cite{sto2016}{\rm)}
Let $\phi\in C(\mathbb{S}^{n-1},\IR^n)$, and let $F$ be defined as follows:
$$F(x)=\begin{cases}
\displaystyle P_{h}[\phi](x),\;\;x\in\mathbb{B}^{n}, \\
\displaystyle \;\;\;\phi(x),\;\;\;\; \;x\in\mathbb{S}^{n-1}.
\end{cases}%\ee
$$
Then
$(1)$
$F$ is hyperbolically harmonic in $\mathbb{B}^{n}$ and continuous in $\overline{\mathbb{B}^{n}}$;\smallskip

\noindent $(2)$
 $\|F\|_{\infty}=\|\phi\|_{\infty}$, where $\|F\|_{\infty}=\sup\{|F(x)|:\; x\in \overline{\IB^n}\}$ and $\|\phi\|_{\infty}=\sup\{|\phi(\xi)|:\; \xi\in \IS^{n-1}\}$.

Conversely, if $H$ is hyperbolically harmonic in $\mathbb{B}^{n}$ and continuous in $\overline{\mathbb{B}^{n}}$, then $$H=P_{h}[H].$$
\end{Thm}

%%%%%%%%%%%%%%%%%%%%%%%%%%%%%
\section{Representation of solutions to $\Delta_hu=\psi$}\label{sec-2}
%%%%%%%%%%%%%%%%%%%%%%%%%%%%%

The main purpose of this section is to prove Theorem \ref{thm-3}. In this section, we always assume that $n\geq 3$. Before the proof, we recall the following results.

\begin{Thm}\label{xjl-12} {\rm(}\cite[Lemma 3.2]{sto2012}{\rm)} If $u\in C^2(\mathbb{B}^{n})$, then

%\be\label{eq03.02.0}
$$u(0)=\int_{\mathbb{S}^{n-1}}u(r\xi)\,d\sigma(\xi)-\int_{\mathbb{B}^{n}(0,r)}g(|x|,r)\Delta_{h}u(x)\,d\tau(x),$$
where $0<r<1$ and $g(t,r)$ is defined in \eqref{eq0.0.01}.
\end{Thm}

\begin{Thm}\label{jl-2} {\rm(}\cite[Corollary 4.1]{sto2012}{\rm)}
 For any $y\in \IB^n$, $$\int_{\mathbb{B}^{n}}G_{h}(x,y)d\nu(x)=\frac{1}{2n(n-1)}(1-|y|^{2})^{n-1}.$$
\end{Thm}

The next two theorems are about the M\"{o}bius invariance of $P_{h}[f]$ and $d\tau $.

\begin{Thm}\label{Thm-C} {\rm(}\cite[Theorem 5.23]{sto1999} {\rm or} \cite{sto2016}{\rm)}
If $f\in L^{1}(\mathbb{S}^{n-1})$, then $$P_{h}[f\circ \varphi]=P_{h}[f]\circ \varphi$$ for all $\varphi\in \mathcal{M}(\mathbb{B}^{n})$.
\end{Thm}

\begin{Thm}\label{Thm-M} {\rm(}\cite[Theorem 3.4(b)]{sto1999} {\rm or} \cite{sto2016}{\rm)}
If $f\in L^{1}(\mathbb{B}^{n},\tau)$ and $\varphi\in \mathcal{M}(\mathbb{B}^{n})$, then
 $$ \int_{\mathbb{B}^{n}}f(x) d\tau(x) =\int_{\mathbb{B}^{n}}f\circ \varphi (y) d\tau(y) .$$
\end{Thm}

\begin{lem}\label{lem1-10}
Let $\psi\in C(\mathbb{B}^{n},\IR^n)$. Suppose there is a constant $\mu_{1}$ such that $$\int_{\mathbb{B}^{n}}(1-|x|^{2})^{n-1} |\psi(x)|\,d\tau(x)\leq \mu_{1}.$$ Then $$\int_{\mathbb{B}^{n}}g(|x|) |\psi(x)|\,d\tau(x)\leq \mu_3,$$ where $g(r)$ is defined in \eqref{eq0.0.01},
 $\mu_3=\mu_3(n, \mu_{1}, \|\psi\|_{\frac{1}{2},\infty})$ and $\|\psi\|_{\frac{1}{2},\infty} =\sup \{|\psi(x)|:x\in \mathbb{B}^{n}(0,\frac{1}{2})\}$.
\end{lem}
\bpf
By letting $y=0$ in \eqref{eq0.0.210} and Theorem \Ref{jl-2}, we get
$$\int_{\mathbb{B}^{n}(0,\frac{1}{2})}g(|x|)\,d\tau(x)=\int_{\mathbb{B}^{n}(0,\frac{1}{2})}G_{h}(x,0)\frac{d\nu(x)}{(1-|x|^{2})^{n}}\leq\frac{1}{2n(n-1)(1 -\frac{1}{4})^{n}},$$ and thus
the assumption ``$\psi\in C(\mathbb{B}^{n},\IR^n)$" gives that $$\int_{\mathbb{B}^{n}(0,\frac{1}{2})}g(|x|)  |\psi(x)|\,d\tau(x)\leq  \frac{\|\psi\|_{\frac{1}{2},\infty} }{2n(n-1)(1 -\frac{1}{4})^{n}}.$$

Obviously, for $x\in \mathbb{B}^{n}\setminus \mathbb{B}^{n}(0,\frac{1}{2})$,
$$g(|x|)=\frac{1}{2n}\int_{|x|}^{1}\frac{(1-s^{2})^{n-2}}{s^{n}}ds^{2}\leq
\frac{2^{n-1}}{n} \int_{|x|}^{1} (1-s^{2})^{n-2}ds^{2} =
\frac{2^{n-1}}{n(n-1)}(1-|x|^{2})^{n-1},$$ it follows that
$$\int_{\mathbb{B}^{n}\setminus \mathbb{B}^{n}(0,\frac{1}{2})}g(|x|) |\psi(x)|\,d\tau(x)\leq\frac{2^{n-1}}{n(n-1)}\mu_{1}.$$

Since
$$\int_{\mathbb{B}^{n}}g(|x|) |\psi(x)|\,d\tau(x)=\int_{\mathbb{B}^{n}(0,\frac{1}{2})}g(|x|)|\psi(x)|\,d\tau(x)+\int_{\mathbb{B}^{n}\setminus \mathbb{B}^{n}(0,\frac{1}{2})}g(|x|) |\psi(x)|\,d\tau(x),$$
by letting $$\mu_3=\frac{\|\psi\|_{\frac{1}{2},\infty} }{2n(n-1)(1 -\frac{1}{4})^{n}}+ \frac{2^{n-1}}{n(n-1)}\mu_{1},$$ we see that
the lemma holds.
\epf

\begin{lem}\label{lem2-10}
Suppose that $u\in C^{2}(\mathbb{B}^{n},\mathbb{R}^{n}) \cap C(\overline{\mathbb{B}^{n}},\IR^n)$ and satisfies \eqref{eq-1.1}. If $$\int_{\mathbb{B}^{n}}(1-|x|^{2})^{n-1}|\psi(x)|\,d\tau(x)\leq \mu_{1},$$ then
$$u(0)=P_h[\phi](0)-\int_{\mathbb{B}^{n}}G_{h}(0,x)\psi(x)\,d\tau(x).$$
\end{lem}
\bpf It follows from the assumption ``$ \int_{\mathbb{B}^{n}}(1-|x|^{2})^{n-1} |\psi(x)|\,d\tau(x)\leq \mu_{1}$" and Lemma \ref{lem1-10} that $$\int_{\mathbb{B}^{n}}g(|x|)  |\psi(x)|\,d\tau(x)\leq \mu_3.$$
Since
 $$\int_{\mathbb{B}^{n}(0,r)}g(|x|,r)|\psi(x)|\,d\tau(x)\leq \int_{\mathbb{B}^{n}}g(|x| )|\psi(x)|\,d\tau(x),$$
by Lebesgue's Dominated Convergence Theorem, we have that
\be\label{eq03.03.0}\lim_{r\rightarrow 1^{-}}\int_{\mathbb{B}^{n}(0,r)}g(|x|,r) \psi(x) \,d\tau(x)=\int_{\mathbb{B}^{n}}g(|x| ) \psi(x) \,d\tau(x).\ee
Furthermore, the assumption ``$u\in C(\overline{\mathbb{B}^{n}},\IR^n )$" gives that
\be \label{eq10.1.0}
\lim_{r\rightarrow 1^{-}}\int_{\mathbb{S}^{n-1}}u(r\xi)\,d\sigma(\xi)=\int_{\mathbb{S}^{n-1}}u(\xi)\,d\sigma(\xi).\ee
Then Theorem \Ref{xjl-12}, \eqref{eq03.03.0} and \eqref{eq10.1.0} imply that
\beq\nonumber
u(0)&=&\int_{\mathbb{S}^{n-1}}u(\xi)\,d\sigma(\xi)-\int_{\mathbb{B}^{n}}g(|x|)\Delta_{h}u(x)\,d\tau(x) \\ \nonumber
&=&\int_{\mathbb{S}^{n-1}}P_{h}(0,\xi)\phi(\xi)\,d\sigma(\xi)-\int_{\mathbb{B}^{n}}G_h(0,x)\psi(x)\,d\tau(x),
\eeq
 as required.
\epf

\subsection*{Proof of Theorem \ref{thm-3}}
We prove this theorem by two steps.
In the first step, we check that for any fixed $\zeta\in \overline{\mathbb{B}^{n}}(0,r_{0})$, $u\circ \varphi_{\zeta}$ satisfies the requirements in Lemma \ref{lem2-10},
 where $0\leq r_{0}<1$.
 In the second step, by applying Lemma \ref{lem2-10} to $u\circ \varphi_{\zeta}$, we finish the proof.
 \bcl\label{jl0-3} For any fixed $\zeta\in \overline{\mathbb{B}^{n}}(0,r_{0})$, $u\circ \varphi_{\zeta}$ satisfies the requirements in Lemma \ref{lem2-10}. \ecl

Obviously, for any fixed $\zeta\in \overline{\mathbb{B}^{n}}(0,r_{0})$,  $u\circ \varphi_{\zeta} \in C^{2}(\mathbb{B}^{n},\mathbb{R}^{n}) \cap C(\overline{\mathbb{B}^{n}},\IR^n )$. The M\"obius invariance property \eqref{eq0.0.1} and the assumption ``$u\in C(\overline{\mathbb{B}^n},\IR^n)$"
imply that $$\Delta_{h}(u\circ \varphi_{\zeta})(y)=\Delta_{h}u\big( \varphi_{\zeta}(y)\big)=\psi(\varphi_{\zeta}(y))=\psi\circ \varphi_{\zeta}(y)$$ in $\mathbb{B}^{n}$ and
$$(u\circ \varphi_{\zeta})\mid_{\mathbb{S}^{n-1}}=\phi\circ \varphi_{\zeta}.$$
So, to prove the claim, it suffices to show the following: There exists a constant $\mu_{4}$ such that
\be\label{eq5.1.1}\int_{\mathbb{B}^{n}}(1-|y|^{2})^{n-1}|\psi\big(\varphi_{\zeta}(y)\big)|\,d\tau(y)\leq \mu_{4},\ee
where $\mu_{4}=\mu_{4}(\mu_{1},n,r_{0})$.

Let $w=\varphi_{\zeta}(y)$. Then we have that $y=\varphi_{\zeta}(w)$, so Theorem \Ref{Thm-M}, \eqref{wed-1} and \eqref{eq1.0.3}, together with the assumption ``$\int_{\mathbb{B}^{n}}(1-|x|^{2})^{n-1}|\psi(x)|\,d\tau(x)\leq \mu_{1}$", yield
 \begin{eqnarray*}
   &\;\;& \int_{\mathbb{B}^{n}}(1-|y|^{2})^{n-1}|\psi\big(\varphi_{\zeta}(y)\big)|\,d\tau(y)\\ \nonumber
  &=& \int_{\mathbb{B}^{n}}(1-|\varphi_{\zeta}(w)|^{2})^{n-1}|\psi(w)|\,d\tau(w) \;\;\;\;\;\;\;\;\;\;\text{(letting\;}w=\varphi_{\zeta}(y)\text{)}\\ %\nonumber
  &\leq&2^{n-1}\int_{\mathbb{B}^{n}}\frac{ (1-|w|^{2})^{n-1}|\psi(w)|}{(1-|\zeta|)^{n-1}}\,d\tau(w)\;\;\;\;\;\;\;\;\text{(by \eqref{wed-1} and \eqref{eq1.0.3})}\\ %\nonumber
  &\leq& \frac{2^{n-1}\mu_1}{(1-r_{0})^{n-1}}.
 \end{eqnarray*}
 Obviously, letting $\mu_4 = 2^{n-1}\mu_1(1-r_{0})^{1-n}$ yields \eqref{eq5.1.1}.

\bcl\label{jl00-3}
$u =P_{h}[\phi] -G_{h}[\psi]$.
\ecl

By replacing $u$ with $u\circ \varphi_{\zeta}$ and by using \eqref{eq0.0.1} and Theorem \Ref{Thm-C}, we see from
 Lemma \ref{lem2-10} that
  \begin{eqnarray*}
   u(\zeta)&=&u\circ\varphi_{\zeta}(0)=P_{h}[\phi\circ \varphi_{\zeta}](0)-\int_{\mathbb{B}^{n}}G_h(0,y)\Delta_{h}(u\circ \varphi_{\zeta})(y)\,d\tau(y)\\
&=&P_{h}[\phi]\big( \varphi_{\zeta}(0)\big)-\int_{\mathbb{B}^{n}}G_h(0,y)\Delta_{h}u\big(\varphi_{\zeta}(y)\big)\,d\tau(y).
 \end{eqnarray*}
 Let $w=\varphi_{\zeta}(y) $. It follows from $$\varphi_{\zeta}(0)=\zeta, \;\; G_h(0,\varphi_{\zeta}(w))=g(|\varphi_{\zeta}(w)|)=G_h(\zeta,w)$$ and Theorem \Ref{Thm-M} that
 $$u(\zeta)=P_{h}[\phi](\zeta)-\int_{\mathbb{B}^{n}}G_h(\zeta,w)\Delta_{h}u(w)\,d\tau(w)=P_{h}[\phi](\zeta)-G_{h}[\psi](\zeta).$$
By the arbitrariness of $r_0$ in $[0, 1)$, we see that the proof of the theorem is complete.\qed

\section{Lipschitz continuity of $\Phi=P_{h}[\phi]$}\label{LipP}

The aim of this section is to prove the $\omega$-Lipschitz continuity of $\Phi=P_{h}[\phi]$ (Theorem \ref{thm-11}).

Before the proof of Theorem \ref{thm-11}, we need an estimate on $\|D \Phi(x) \|$ in terms of $\omega(1-|x|)$ which is formulated in Lemma \ref{jl23-3}.
The proof of Lemma \ref{jl23-3} needs some preparation which consists of three lemmas. The first lemma is as follows.

 \begin{lem}\label{lem-22}  Suppose $\phi \in C(\mathbb{S}^{n-1},\IR^n)$. Then for each $k\in \{1,2,\ldots,n\}$,\ben
 \item
$\frac{\partial}{\partial x_{k}} \Phi(x)$ is continuous in $\mathbb{B}^{n}$;
\item
$\frac{\partial}{\partial x_{k}} \Phi(x) =\int_{\mathbb{S}^{n-1}}\frac{\partial}{\partial x_{k}}  P_{h}(x,\xi)\phi(\xi) \,d\sigma(\xi)$ for $x\in \mathbb{B}^{n}$.
\een\end{lem}
\bpf
 In order to prove this lemma, we only need to discuss the case
  $k=1$ since other cases can be discussed in a similar way.
For this, we assume that $x\in \overline{\mathbb{B}^{n}}(0,r_{0})$ and $x+\Delta x_{1}\in \overline{\mathbb{B}^{n}}(0,r_{0})$, where $x=(x_1, \ldots, x_n)$, $x+\Delta x_{1}=(x_1+\Delta x_{1}, \ldots, x_n)$ and $0<r_{0}<1$.
 Then $$\frac{\Phi(x+\Delta x_{1})- \Phi(x)}{\Delta x_{1}}= \int_{\mathbb{S}^{n-1}} \frac{ P_{h}(x+\Delta x_{1},\xi)-P_{h}(x,\xi)}{ \Delta x_{1} }\phi(\xi) \,d\sigma(\xi).$$

 Obviously, $\frac{\partial}{\partial x_{1}}P_{h}(x,\xi)\phi(\xi)$ is continuous
 in $\overline{\mathbb{B}^{n}} (0,r_{0})\times \mathbb{S}^{n-1}$, and so
 $$\int_{\mathbb{S}^{n-1}}\frac{\partial}{\partial x_{1}}P_{h}(x,\xi)\phi(\xi)d\sigma(\xi)$$
 is continuous on $ \overline{\mathbb{B}^{n}}(0,r_{0})$.
 By applying the Lagrange mean-value theorem
 to $P_{h}(x,\xi)$ w.r.t. $x_{1}$, we see that there exists $t_1\in (0,1)$ such that
 \begin{eqnarray*}
  \frac{\partial \Phi(x)}{\partial x_{1}}
  &=& \lim_{\Delta x_{1}\rightarrow 0}  \frac{ \Phi(x+\Delta x_{1})- \Phi(x)}{ \Delta x_{1}}= \lim_{\Delta x_{1}\rightarrow 0}  \int_{\mathbb{S}^{n-1}}  \frac{\partial}{\partial x_{1}} P_{h}(x+t_{1}\Delta x_{1},\xi) \phi(\xi) \,d\sigma(\xi)\\
  &=& \int_{\mathbb{S}^{n-1}}  \frac{\partial}{\partial x_{1}} P_{h}(x,\xi) \phi(\xi) \,d\sigma(\xi),\end{eqnarray*}
 as required.
 \epf

Let $\xi_{0}=e_{1}\in \mathbb{S}^{n-1}$ denote the first unit coordinate vector $(1,0,\ldots,0)$. Then we have the following estimate.
 \begin{lem}\label{cjlem-22} Suppose $q\geq0$, $p-q-n>0$ and $n\geq 3$. Then

  \beq\label{eqcj1.1.3}
 & \;\;&  \int_{\mathbb{S}^{n-1}} \frac{|\xi-\xi_{0}|^{q}\omega(|\xi-\xi_{0}|)}{|\xi -r\xi_{0} |^{p}} \,d\sigma(\xi)\\  \nonumber
& \leq& \alpha_1 \frac{\omega(1-r)}{(1-r)^{p+1-q-n}}\left( \frac{1}{ q+n-1 } + \frac{2^{p}}{p-q-n} \right),
\eeq
where $\omega$ is a majorant, $0\leq r< 1$ and $\alpha_1= \omega_{n-2}/\omega_{n-1} =\frac{1}{\sqrt{\pi}}\frac{\Gamma(\frac{n}{2})}{\Gamma(\frac{n-1}{2})}.$
\end{lem}
\bpf We shall prove this lemma by using a similar argument as in \cite{AKM2008} and \cite{arsen2009}.
In order to estimate the integral in \eqref{eqcj1.1.3}, we split $\mathbb{S}^{n-1}$ into the following two subsets:
$$E=\{\xi\in \mathbb{S}^{n-1}:|\xi-\xi_{0}|\leq 1-r\}\;\; \mbox{and}\;\; F=\{\xi\in \mathbb{S}^{n-1}:|\xi-\xi_{0}|> 1-r\}.$$
Then \eqref{eqcj1.1.3} easily follows from the following two claims.
 \bcl\label{eq1.1.5} $\int_{E} \frac{ |\xi-\xi_{0}|^{q}\omega(|\xi-\xi_{0}|)}{ |\xi- r\xi_{0} |^{p}} \,d\sigma(\xi)\leq \alpha_1\frac{\omega(1-r)}{ q+n-1}(1-r)^{q+n-1-p}.$
\ecl

Since $|\xi-r\xi_{0}|\geq 1-|r\xi_{0}|=1-r$ for all $\xi\in \mathbb{S}^{n-1}$, we have
\beq\label{eq1.1.3}
\int_{E} \frac{|\xi-\xi_{0}|^{q}\omega(|\xi-\xi_{0}|)}{|\xi -r\xi_{0} |^{p}} \,d\sigma(\xi)&\leq& \int_{E} \frac{|\xi-\xi_{0}|^{q}\omega(|\xi-\xi_{0}|)}{(1-r)^{p}} \,d\sigma(\xi) \\  \nonumber
&=&\frac{(1-r)^{-p}}{\omega_{n-1}  }\int_{E}  |\xi-\xi_{0}|^{q} \omega(|\xi-\xi_{0}|)\,dS(\xi), \eeq
where $S$ denotes the $(n-1)$-dimensional Lebesgue measure on $\IS^{n-1}$. Let $\xi=(\xi_{1},\ldots,\xi_{n})\in E$ has the expression \eqref{Ttrans}. Then, $\theta_{1}\in [0,\varphi_r]\subset [0,\frac{\pi}{2}]$, $\theta_{2},\ldots,\theta_{n-2}\in [0,\pi]$ and $\theta_{n-1}\in [0,2\pi]$, where $\varphi_r=2\arcsin\frac{1-r}{2}$.
It follows from \eqref{eq1.0.0} that
 \begin{eqnarray*}
\int_{E}  |\xi-\xi_{0}|^{q} \omega(|\xi-\xi_{0}|)\,dS(\xi)
 &=&\int_{0}^{\varphi_r} (2-2\cos\theta_{1})^{\frac{q}{2}} \omega\big((2-2\cos\theta_{1})^{\frac{1}{2}}\big ) \sin^{n-2}\theta_{1}\,d\theta_{1}\\
 &\;\;&\cdot\int_{0}^{\pi} \sin^{n-3}\theta_{2}\,d\theta_{2}\cdots \int_{0}^{\pi} \sin\theta_{n-2}\,d\theta_{n-2} \int_{0}^{2\pi} \, d\theta_{n-1}   \\
  &=&\omega_{n-2} \int_{0}^{\varphi_r} (2-2\cos\theta_{1})^{\frac{q}{2}}  \omega\big((2-2\cos\theta_{1})^{\frac{1}{2}} \big) \sin^{n-2}\theta_{1}\,d\theta_{1}.
 \end{eqnarray*}
Let $\rho =\sqrt{2-2\cos\theta_{1}}\in [0,1-r]$. Then $d\theta_{1}=\frac{\rho}{\sin \theta_{1}} d\rho$, from which we deduce that
\beq\label{eq1.1.4}
 & \;\;&\int_{E}  |\xi-\xi_{0}|^{q}  \omega(|\xi-\xi_{0}|)\,dS(\xi)\\  \nonumber
& \leq&
\omega_{n-2} \int_{0}^{1-r} \rho^{q+n-2} \omega(\rho)\,d\rho\leq\frac{\omega_{n-2}\omega(1-r) }{q+n-1}(1-r)^{q+n-1},\eeq where in the first inequality, the relation
$\sin^{2}\theta_{1}=\rho^{2}\left(1-\frac{\rho^{2}}{4}\right)\leq \rho^{2}$ is applied.
It follows from $\alpha_1=\omega_{n-2}/\omega_{n-1} $, \eqref{eq1.1.3} and \eqref{eq1.1.4} that
\[
 \int_{E} \frac{ |\xi-\xi_{0}|^{q}\omega(|\xi-\xi_{0}|)}{ |\xi- r\xi_{0} |^{p}} \,d\sigma(\xi)\leq \alpha_1\frac{\omega(1-r)}{ q+n-1}(1-r)^{q+n-1-p},
\]
which is what we need.
\bcl \label{eq1.1.6}
$\int_{F} \frac{|\xi-\xi_{0}|^{q}}{|r\xi_{0}-\xi |^{p}}\omega(|\xi-\xi_{0}|) \,d\sigma(\xi)\leq 2^{p}\alpha_1\frac{\omega(1-r) }{p-q-n}(1-r)^{q+n-1-p}.$
 \ecl

Since for all $\xi\in F$,
$$|\xi-\xi_{0}|\leq |\xi-r\xi_{0}|+|r\xi_{0}-\xi_{0}|=|\xi-r\xi_{0}|+1-r\;\;\mbox{and}\;\; |\xi-r\xi_{0}|\geq 1-r,$$ we easily see that
$$|\xi-\xi_{0}|\leq 2|\xi-r\xi_{0}|,$$ and so
$$\int_{F} \frac{|\xi-\xi_{0}|^{q}\omega(|\xi-\xi_{0}|)}{|\xi -r\xi_{0}|^{p}} \,d\sigma(\xi)\leq 2^{p} \int_{F}  |\xi-\xi_{0}|^{q-p} \omega(|\xi-\xi_{0}|)\,d\sigma(\xi).$$
Then the similar reasoning as in the proof of \eqref{eq1.1.4} leads to
\begin{eqnarray*}
\int_{F} \frac{|\xi-\xi_{0}|^{q}}{|r\xi_{0}-\xi |^{p}}\omega(|\xi-\xi_{0}|) \,d\sigma(\xi)
& \leq& 2^{p} \alpha_1 \int_{1-r}^{2} \rho^{q-p}\omega(\rho)\rho^{n-2}\,d\rho\\  \nonumber
& \leq&  \frac{2^{p}\alpha_1 \omega(1-r) }{p-q-n}(1-r)^{q+n-1-p},
 \end{eqnarray*}
where, in the last inequality, the assumption that $\frac{\omega(t)}{t}$ is non-increasing, is exploited.  Hence Claim \ref{eq1.1.6} is true. \epf

Based on Lemmas \ref{lem-22} and \ref{cjlem-22}, we have the following estimate on $\left |\frac{\partial \Phi}{\partial x_{k}}(x)\right|$.
 \begin{lem}\label{jl10-3} Suppose $\phi$ and $\omega$ satisfy the conditions in Theorem \ref{thm-11}.
Then there is a constant $\alpha_{2}$ such that for all $x\in [0, e_1)$ and $k\in \{1,2,\ldots,n\}$,
$$\left |\frac{\partial \Phi}{\partial x_{k}}(x)\right|\leq \alpha_{2}\frac{\omega(1-r)}{1-r},$$ where $\alpha_{2}=\alpha_{2}(n)$, $[0, e_1)=\{x\in \mathbb{B}^{n}:\; x=re_1, 0\leq r<1\}$ and $n\geq 3$.
\end{lem}
\bpf
For any $x_0\in [0, e_1)$, obviously, there is an $r\in [0, 1)$ such that $x_{0}=r\xi_{0},$ where $\xi_{0}=e_{1}$.
We prove the claim by considering two cases.

 \begin{case}\label{jl12-3} $2\leq k\leq n$. \end{case}
Since \eqref{eq1.1.0} implies
%\beq \label{eq1.1.1}
$$ \frac{\partial}{\partial x_{k}}P_{h}(x_{0},\xi)
  =\frac{2(n-1)(1-|x_{0}|^{2})^{n-1}\xi_{k}}{|\xi-x_{0}|^{2n }},$$
%\eeq
 we infer from
\be \label{eqc1.1.1}\nonumber \phi(\xi_{0})=\int_{\mathbb{S}^{n-1}}P_{h} (x,\xi)\phi(\xi_{0})\,d\sigma(\xi),\ee
 together with Lemma \ref{lem-22}, that
 \begin{eqnarray*}
 \left|\frac{\partial \Phi}{\partial x_{k}}(x_{0})\right|
 &=& \left| \int_{\mathbb{S}^{n-1}} \frac{\partial}{\partial x_{k}}P_{h}(x_{0},\xi) \phi(\xi) \,d\sigma(\xi)\right|\\
  &=& \left|\int_{\mathbb{S}^{n-1}}\frac{2(n-1)(1-|x_{0}|^{2})^{n-1}\xi_{k}}{|\xi-x_{0}|^{2n }}\big(\phi(\xi)-\phi(\xi_{0})\big) \,d\sigma(\xi)\right|\\
 &\leq& 2(n-1)(1-|x_{0}|^{2})^{n-1} \int_{\mathbb{S}^{n-1}} \frac{|\xi_{k}| \cdot \big|\phi(\xi)-\phi(\xi_{0})\big|}{|\xi-x_{0}|^{2n }} \,d\sigma(\xi) .
 \end{eqnarray*}
By using the fact $|\xi_{k}|\leq |\xi-\xi_{0}|$ for $2\leq k\leq n$ and the assumption ``$|\phi(\xi)-\phi(\xi_{0})|\leq \omega(|\xi-\xi_{0}|)$", we get
%\be\label{eq1.1.2}
$$ \left|\frac{\partial \Phi}{\partial x_{k}}(x_{0})\right|\leq 2(n-1)(1-|x_{0}|^{2})^{n-1} \int_{\mathbb{S}^{n-1}}\frac{| \xi-\xi_{0}| \omega(|\xi-\xi_{0}|)}{|\xi-x_{0}|^{2n } }\,d\sigma(\xi) .
 $$
Then Lemma \ref{cjlem-22} leads to
$$\left|\frac{\partial \Phi}{\partial x_{k}}(re_{1})\right| \leq \alpha_{3}\frac{\omega(1-r)}{1-r},
$$
where $\alpha_3=\big(\frac{ n-1}{n}+2^{2n}\big)2^{n}\alpha_1$.

 \begin{case}\label{jl13-3} $k=1$.\end{case}
Again, \eqref{eq1.1.0} implies
$$
\frac{\partial}{\partial x_{1}}P_{h}(re_{1},\xi)=\frac{2(n-1)(1-|x_{0}|^{2})^{n-1}(\xi_{1}-|x_{0}|)}{|x_{0}-\xi|^{2n}}-\frac{2(n-1)(1-|x_{0}|^{2})^{n-2}|x_{0}|}{|x_{0}-\xi|^{2n-2}},
$$
and so
 \begin{eqnarray*}
 \left|\frac{\partial \Phi}{\partial x_{1}}(re_{1})\right|
 &=& \left|\int_{\mathbb{S}^{n-1}} \frac{\partial}{\partial x_{1}}P_{h}(x_{0},\xi)\big(\phi(\xi)-\phi(\xi_{0})\big) \,d\sigma(\xi)\right|\\
 &\leq&  2(n-1)(1-|x_{0}|^{2})^{n-1}\int_{\mathbb{S}^{n-1}} \frac{\omega(|\xi-\xi_{0}|)
 \big|\xi_{1}-|x_{0}|\big| }{|x_{0}-\xi|^{2n }} \,d\sigma(\xi)\\
  &\;\;& + 2(n-1)(1-|x_{0}|^{2})^{n-2}|x_{0}|\int_{\mathbb{S}^{n-1}} \frac{\omega(|\xi-\xi_{0}|)}{|x_{0}-\xi|^{2n-2 }} \,d\sigma(\xi).
 \end{eqnarray*}
Since $$\big|\xi_{1}-|x_{0}|\big|\leq  |\xi_{1}-1|+\big|1-|x_{0}|\big|\leq |\xi-\xi_{0}|+1-|x_{0}|,$$
 we get
%\beq \label{eq1.1.7}\;\;\;\;
 \begin{eqnarray*}
 \left|\frac{\partial \Phi}{\partial x_{1}}(re_{1})
 \right|
 &\leq&  2(n-1)(1-|x_{0}|^{2})^{n-1}\int_{\mathbb{S}^{n-1}} \frac{\omega(|\xi-\xi_{0}|)|\xi-\xi_{0}|}{|x_{0}-\xi|^{2n }} \,d\sigma(\xi)\\ %\nonumber
    &\;\;&+ 2(n-1)(1-|x_{0}|)^{n}(1+|x_{0}|)^{n-1} \int_{\mathbb{S}^{n-1}} \frac{\omega(|\xi-\xi_{0}|)}{|x_{0}-\xi|^{2n }} \,d\sigma(\xi)\\ %\nonumber
&\;\;&+2(n-1)(1-|x_{0}|^{2})^{n-2}|x_{0}|\int_{\mathbb{S}^{n-1}} \frac{\omega(|\xi-\xi_{0}|) }{|x_{0}-\xi|^{2n-2 }} \,d\sigma(\xi).
  \end{eqnarray*}%\eeq
Furthermore, Lemma \ref{cjlem-22} guarantees that
$$
\int_{\mathbb{S}^{n-1}} \frac{ \omega(|\xi-\xi_{0}|)|\xi-\xi_{0}|}{|x_{0}-\xi|^{2n}} \,d\sigma(\xi)
\leq \alpha_1\frac{\omega(1-r)}{(1-r)^{n}}\left(\frac{1}{n}+\frac{4^{n}}{n-1}\right),
$$
$$
\int_{\mathbb{S}^{n-1}} \frac{\omega( |\xi-\xi_{0}| )}{|x_{0}-\xi|^{2n}} \,d\sigma(\xi)
\leq \alpha_1\frac{\omega(1-r)}{(1-r)^{n+1}}\left(\frac{1}{n-1}+\frac{4^{n}}{n}\right)
$$ and
$$
\int_{\mathbb{S}^{n-1}} \frac{\omega( |\xi-\xi_{0}| ) }{|x_{0}-\xi|^{2n-2 }} \,d\sigma(\xi)
\leq \alpha_1\frac{\omega(1-r)}{(1-r)^{n-1}}\left(\frac{1}{n-1}+\frac{4^{n-1}}{n-2}\right).
$$
 Hence $$ \left|\frac{\partial \Phi}{\partial x_{1}}(re_{1})\right|\leq \alpha_{4}\frac{\omega(1-r)}{1-r},
$$
where
\[
\alpha_4=(n-1)\left( \frac{2^{n}+8^{n}}{n}+\frac{2^{n }+8^{n}+2^{n-1}}{n-1 } +\frac{8^{n-1}}{ n-2}\right)\alpha_{1}.
\]\medskip

By letting $\alpha_2=\max\left\{\alpha_3, \alpha_4\right\}$, we see that the lemma is true.
 \epf

Now, we are ready to state and prove the main lemma in this section.
\begin{lem}\label{jl23-3}  Suppose $\phi$ and $\omega$ satisfy the conditions in Theorem \ref{thm-11}.
Let $\alpha_{0}=\sqrt{n}\alpha_{2}$, where $n\geq 3$ and $\alpha_{2}=\alpha_{2}(n)$ is the same constant as in Lemma \ref{jl10-3}. Then
$$\|D \Phi(x) \|\leq \alpha_{0} \frac{\omega(1-|x|)}{1-|x|}$$ in $\mathbb{B}^{n}$.
\end{lem}
\bpf Let $x_0\in \mathbb{B}^n$.
We divide the proof into two cases.

\bca\label{wed-2} $x_{0}\in [0,e_{1})$.\eca
Since the Cauchy-Schwarz inequality implies
\begin{eqnarray*}
\|D \Phi(x_{0})\|&=&\sup_{\zeta\in \mathbb{S}^{n-1}}|D \Phi(x_{0}) \zeta|
=\sup_{\zeta\in \mathbb{S}^{n-1} }\left| \sum_{k=1}^{n}\frac{\partial \Phi}{\partial x_{k}}(x_{0}) \cdot \zeta_{k} \right| \leq  \left(\sum_{k=1}^{n}\left|\frac{\partial \Phi}{\partial x_{k}}(x_{0}) \right|^{2}\right)^{\frac{1}{2}},
 \end{eqnarray*}
we see from Lemma \ref{jl10-3} that
\be\label{eqc1.c1.c1}\|D \Phi(x_{0})\|\leq\left(\sum_{k=1}^{n}\left|\frac{\partial \Phi}{\partial x_{k}}(x_{0}) \right|^{2}\right)^{\frac{1}{2}}   \leq \sqrt{n}\alpha_{2} \frac{\omega(1-|x_{0}|)}{1-|x_{0}|}.\ee
 The proof of Lemma \ref{jl23-3} holds in this case.

\bca $x_{0}\notin [0,e_{1})$.
\eca
For the proof in this case, we choose a unitary transformation $U$ such that $U(re_{1})=x_{0}$, $r=|x_0|$, and for $y\in\mathbb{B}^{n}$, let $$W(y)=:\Phi(U(y)).$$ By Theorem \Ref{Thm-C}, we see that $$W=P_h[\phi]\circ U=P_h[\phi\circ U].$$
 Then we have the following claim.
\bcl\label{c-1} $\|D W(re_{1})\| \leq \sqrt{n}\alpha_{2} \frac{\omega(1-r)}{1-r}.$\ecl

%We prove this Claim by two steps. First, we will prove that $\phi\circ U$ satisfies the condition in Theorem \ref{thm-11}.
%Then, \eqref{eqc1.c1.c1} implies the result.
The assumption ``$|\phi(\xi)-\phi(\eta)|\leq \omega(|\xi-\eta|)$" implies that for $\xi$, $\eta\in \mathbb{S}^{n-1}$,
$$|\phi(U(\xi))-\phi(U(\eta))|\leq \omega(| U(\xi) - U(\eta) |)=\omega(|\xi-\eta|).$$
Thus, by replacing $\Phi$ by $\Phi\circ U$, the similar reasoning as in the discussions of Case \ref{wed-2} shows that
$$\|D W(re_{1})\| \leq \sqrt{n}\alpha_{2} \frac{\omega(1-r)}{1-r},$$
which is what we want.\medskip

Now, we are ready to finish the proof of the lemma in this case.
By applying the chain rule, we obtain
 $$D W(y)\big|_{ y=re_{1}} =D (\Phi\circ U)(y)\big|_{y= re_{1}}=
  (D \Phi) \circ U(y)\big|_{y= re_{1}} \times D U(y) \big|_{y= re_{1}}
 =D \Phi( x_{0}) \times U,$$
where $\times$ denotes the usual matrix product. Then
\beq \label{eq15.15.17}
 &\;\;&\|D\Phi( x_{0})\|=\|D\Phi( x_{0}) \times  U\|=\|D W(re_{1})\|
    \leq\sqrt{n}\alpha_{2}(n) \frac{\omega(1-r)}{1-r}.
 \eeq

 By \eqref{eqc1.c1.c1} and \eqref{eq15.15.17}, we complete the proof of  Lemma \ref{jl23-3}.
\epf

\subsection*{Proof of Theorem \ref{thm-11}}
Now, we are ready to prove Theorem \ref{thm-11} by applying Lemma \ref{jl23-3}. \medskip

Lemma \ref{lem-22}(a) implies that $\Phi\in C^{1}(\mathbb{B}^{n})$ and so $\Phi$ is differentiable.
Since $\mathbb{B}^{n}$ is a $\Lambda_{w}$-extension domain for a fast majorant $\omega$ (cf. \cite[Section 1]{dyak2004}), it follows from
the mean-value theorem of differentials  (see e.g. \cite[Theorem 9.19]{wrudin}),
\eqref{eq21.21.221} and Lemma \ref{jl23-3} that
there is a rectifiable curve $\gamma \subset \mathbb{B}^{n}$ joining $x$ to $y$ satisfying
\begin{eqnarray*}
|\Phi(y)-\Phi(x)|\leq \int_{\gamma}\left\|D \Phi(\zeta)\right\|\, ds(\zeta)
\leq \alpha_{0}\int_{\gamma} \frac{\omega(\delta_{\mathbb{B}^{n}}(\zeta))}{\delta_{\mathbb{B}^{n}}(\zeta)} ds(\zeta)
\leq C\ \alpha_{0}\ \omega(|x-y|),
 \end{eqnarray*}
since $\delta_{\mathbb{B}^{n}}(\zeta)=1-|\zeta|$ for $\zeta\in \mathbb{B}^{n}$, where $C=C(\mathbb{B}^{n},\omega)$
 is the same constant as in \eqref{eq21.21.221}. So the proof of this Theorem \ref{thm-11} is complete. \qed

\section{Lipschitz continuity of $\Psi =G_{ h}[\psi]$}\label{LipG}

%\subsubsection{Auxiliary results}
In this section, Theorem \ref{thm-2} is proved through a series of lemmas. From this and Theorem \ref{thm-11}, we derive Lipschitz continuity of $u=P_{h}[\phi]-G_{ h}[\psi]$, i.e. Theorem \ref{thm-4}.

First, let us recall the following lemma from \cite{reka}.

\begin{Thm}\label{Lem-D} {\rm(}\cite[Section 2]{reka}{\rm)}
Let $f$ be a continuous function in $[-1,1]$. Then for any $\eta \in \mathbb{S}^{n-1}$ and $n\geq 3$,
$$\int_{\mathbb{S}^{n-1}}f(\langle \xi,\eta\rangle)\,d\sigma(\xi)=\frac{\Gamma(\frac{n}{2})}{\Gamma(\frac{n-1}{2})\Gamma(\frac{1}{2})} \int_{-1}^{1}(1-t^{2})^{\frac{n-3}{2}}f(t)\,dt.$$
\end{Thm}

\begin{lem}\label{lem-4.4}
Let $$I_{1}(s)=: \int_{0}^{1} F\left(1,\frac{4-n}{2};\frac{n}{2};ts \right )\, dt.$$
\ben
\item
If $n\geq4$, then for $s\in[0,1)$,$$|I_{1}(s)|\leq \mu_{2,2},$$ where $\mu_{2,2}=\mu_2(n,1,4,1)$ is defined in Lemma \ref{lem-4.5};
\item
If $n=3$ and $ s_{0}\in (0,1)$, then for $s\in[0,s_{0}]$, $$|I_{1}(s)|\leq \frac{1}{1-s_{0}}.$$
\een
\end{lem}
\bpf It follows from Lemma \ref{jl-30} that $$I_{1}(s)=\sum_{k=0}^{\infty}\frac{(1)_{k} (\frac{4-n}{2})_{k}}{(\frac{n}{2})_{k}k!}\frac{s^k}{k+1}.$$

If $n\geq4$, the result follows from Lemma \ref{lem-4.5}.

If $n=3$, then for any $s_{0}\in (0,1)$,
we have $$I_{1}(s)=\sum_{k=0}^{\infty}\frac{(1)_{k} (\frac{4-n}{2})_{k}}{(\frac{n}{2})_{k}k!}\frac{s^k}{k+1}\leq \sum_{k=0}^{\infty} s^k \leq \frac{1}{1-s_{0}}.$$
Hence Lemma \ref{lem-4.4} is proved.
\epf

By Lemma \ref{lem-4.4}, we have the following estimate.
\begin{lem}\label{lem1-100}
If $n\geq 3$, then for all $x\in \mathbb{B}^{n} $,
%\be\label{eq2.02.004}
$$J_{n}(x)=:\int_{\mathbb{B}^{n} }\frac{ d\nu(y)}{|y|^{n-2}[x,y]^{2}} \leq \mu_5,$$
where $\mu_5=\max\left\{\frac{ n}{2}  \mu_{2,2}, 65\frac{7}{8}\right\}$, where $\mu_{2,2}$ is the same constant as in Lemma \ref{lem-4.4}.
\end{lem}
\bpf
For $n\geq 3$, \eqref{eq0.0cj.0} leads to
$$
 J_{n}(x)=\int_{\mathbb{B}^{n} }\frac{ d\nu(y)}{|y|^{n-2}[x,y]^{2}}
=n  \int_{0 }^{1}\rho \,d\rho \int_{\mathbb{S}^{n-1} }\frac{\,d\sigma(\xi)}{[x,\rho\xi]^{2}}.
$$
By \eqref{eq1.0.2}, we have
$$\int_{\mathbb{S}^{n-1} }\frac{\,d\sigma(\xi)}{[x,\rho\xi]^{2}}=\int_{\mathbb{S}^{n-1} } \left(1+\rho^{2}|x|^{2}-2\rho|x|\left\langle\frac{x}{|x|},\xi\right\rangle\right)^{-1}\,d\sigma(\xi),$$
and so Theorem \Ref{Lem-D} and \eqref{eq2.2.4} lead to
\beq\nonumber
\int_{\mathbb{S}^{n-1} }\frac{\,d\sigma(\xi)}{[x,\rho\xi]^{2}}
&=& \frac{\Gamma(\frac{n}{2})}{\Gamma(\frac{n-1}{2})\Gamma(\frac{1}{2})}\int_{-1}^{1} (1-t^{2})^{\frac{n-3}{2}}\left(1+\rho^{2}|x|^{2}-2\rho|x|t\right)^{-1}dt\\\nonumber
 &=&F\left(1,\frac{4-n}{2};\frac{n}{2};\rho^{2}|x|^{2}\right).
\eeq
Hence we have
\be\label{eq1.03.30} J_{n}(x) = n \int_{0 }^{1}\rho   F\left(1,\frac{4-n}{2};\frac{n}{2};\rho^{2}|x|^{2}\right)\,d\rho.\ee

When $n\geq 4$, it follows from Lemma \ref{lem-4.4} that for all $x\in \mathbb{B}^{n}$,
\be\label{eq1.cj03.cj30}J_{n}(x) \leq\frac{  n}{2}  \mu_{2,2}.\ee

In the following, we assume that $n=3$. Then we have the following assertion.
 \bcl\label{jl-3} $J_{3}(x)\leq 65\frac{7}{8}.$\ecl

We divide the proof into two cases according to the value of $|x|$.

\bca  $\frac{3}{4}\leq |x|<1$.\eca
 Since $n=3$, by \eqref{eq1.0.2}, we see that
$$ J_{3}(x)=\int_{\mathbb{B}^{3} }\frac{ d\nu(y)}{|y|[x,y]^{2}}\leq  \int_{\mathbb{B}^{3} }\frac{ d\nu(y)}{|y|\cdot |y-x|^{2}}.$$
Let $\delta_{1}= |x|/3$. By \eqref{eq0.0cj.0} and elementary calculations, we have that
$$
\int_{\mathbb{B}^{3}(0,\delta_{1})} \frac{ d\nu(y)}{|y|\cdot |y-x|^{2}}
\leq \int_{\mathbb{B}^{3}(0,\delta_{1})} \frac{ d\nu(y)}{4 |y|\delta_{1}^{2}}=\frac{3}{8},
$$

$$
 \int_{\mathbb{B}^{3}\cap\mathbb{B}^{3}(x,\delta_{1})} \frac{ d\nu(y)}{|y|\cdot |y-x|^{2}}
\leq  \int_{\mathbb{B}^{3}(x,\delta_{1})} \frac{ d\nu(y)}{ 2|y-x|^{2}\delta_{1}}
= \int_{\mathbb{B}^{3}(0,\delta_{1})} \frac{ d\nu(y)}{ 2|y|^{2}\delta_{1} }
=\frac{3}{2}
$$
 and
$$
\int_{\mathbb{B}^{3} \setminus \big( \mathbb{B}^{3}(0,\delta_{1})  \cup \mathbb{B}^{3}(x,\delta_{1})\big)}  \frac{ d\nu(y)}{|y|\cdot |y-x|^{2}}\leq  \int_{\mathbb{B}^{3} \setminus\big( \mathbb{B}^{3}(0,\delta_{1})  \cup \mathbb{B}^{3}(x,\delta_{1})\big)}    \frac{ d\nu(y)}{ \delta_{1}^{3}}\leq\frac{1}{\delta_{1}^{3}}\leq64.
$$
 These inequalities show that the claim holds since
$$\mathbb{B}^{3}= \mathbb{B}^{3}(0,\delta_{1})\cup\mathbb{B}^{3}(x,\delta_{1})\cup \Big(\mathbb{B}^{3} \setminus \big( \mathbb{B}^{3}(0,\delta_{1})  \cup \mathbb{B}^{3}(x,\delta_{1})\big) \Big). $$

\bca $|x|<\frac{3}{4}$.\eca

Under this assumption, we see from \eqref{eq1.03.30} and Lemma \ref{lem-4.4} that
\be\nonumber \label{eq4.2.402}J_{3}(x)= \frac{3 }{2}  \int_{0 }^{1} F\left(1,\frac{1}{2};\frac{3}{2};\rho^{2}|x|^{2}\right)d(\rho^{2})\leq \frac{3 }{2}\cdot \frac{1}{1-\frac{9}{16}}=\frac{24}{7},
\ee as required. So Claim \ref{jl-3} is proved.\medskip

Now, we obtain from \eqref{eq1.cj03.cj30} and Claim \ref{jl-3} that for all $x\in \mathbb{B}^{n}$,
$$J_n(x)=\int_{\mathbb{B}^{n} }\frac{ d\nu(u)}{|u|^{n-2}[x,u]^{2}}\leq\max\left\{\frac{ n}{2}  \mu_{2,2}, 65\frac{7}{8}\right\},$$
and hence the proof of Lemma \ref{lem1-100} is complete.
\epf

Based on Theorem \Ref{Lem-D}, Lemmas \ref{lem-4.4} and \ref{lem1-100},
we obtain some properties of the two unbounded integrals $G_{h} [\psi](x)$ and
\[
\int_{\mathbb{B}^{n}} \frac{\partial G_{h}}{\partial x_{k}}(x,y) \psi(y) \,d\tau (y),
\] which will be presented in the next four lemmas. The first two lemmas deal with the uniform convergence of these two integrals, respectively.
%$G_{h} [\psi](x)$ and $\int_{\mathbb{B}^{n}} \frac{\partial G_{h}}{\partial x_{k}}(x,y) \psi(y) \,d\tau (y)$,

\begin{lem}\label{lem-4.23} Suppose $n\geq3$, $\psi\in C (\mathbb{B}^{n},\mathbb{R}^{n})$ and $|\psi(x)|\leq M(1-|x|^{2})$ in $\mathbb{B}^{n}$, where $M$ is a constant.
Then for all $0<r_{0}<1$, the unbounded integral $ G_{h} [\psi](x)$ is uniformly convergent w.r.t. $x$ in $\overline{\mathbb{B}^{n}}(0,r_{0})$.

\end{lem}
 \bpf  By the assumption ``$|\psi(x)|\leq M(1-|x|^2)$", we see from \eqref{eq20.20.20} that
% \be\label{eq4.2.2}
$$ \big|G_{h} [\psi](x)\big|
 \leq \frac{M}{n}\int_{\mathbb{B}^{n}}\left[ \frac{1}{(1-|y|^2)^{n-1}}\int_{|\varphi_{x}(y)|}^{1}\frac{(1-s^2)^{n-2}}{s^{n-1}}ds\right]\, d\nu(y).$$
For $x\in
 \overline{\mathbb{B}^n}(0,r_{0}) $, \eqref{wed-1} leads to $$[x,y]\geq1-|x|\geq 1-r_0.$$
Since
$$\int_{\mathbb{B}^{n}}\left[ \frac{1}{(1-|y|^2)^{n-1}}\int_{|\varphi_{x}(y)|}^{1}\frac{(1-s^2)^{n-2}}{s^{n-1}}ds\right]\, d\nu(y)=\int_{\mathbb{B}^{n}}  \frac{g(|\varphi_{x}(y)|)}{(1-|y|^2)^{n-1}} \, d\nu(y),$$
we see that

 \begin{eqnarray*}
\;\;\;\;\;\;\;\;\;\;\;\;\;\;&\;&\int_{\mathbb{B}^{n}}\left[ \frac{1}{(1-|y|^2)^{n-1}}\int_{|\varphi_{x}(y)|}^{1}\frac{(1-s^2)^{n-2}}{s^{n-1}}ds\right]\, d\nu(y)\\
 &\leq & \frac{1}{n(n-2)}\int_{\mathbb{B}^{n}}\frac{   (1-|\varphi_{x}(y)|^{2})^{n-1}   }{|\varphi_{x}(y)|^{n-2}(1-|y|^{2})^{n-1}}\, d\nu(y)\;\;\;\;\;\;\text{(by Lemma \ref{eq0.0.21})}\\
  &<& \frac{1}{n(n-2)}\int_{\mathbb{B}^{n}}\frac{1 }{ [x,y]^{n}|x-y|^{n-2}}\, d\nu(y)\;\;\;\;\;\;\;\;\;\;\;\;\;\;\;\;\;\;\;\;\;\;\;\;\text{(by \eqref{eq1.0.3})}\\
 &\leq& \frac{1}{n(n-2)(1-r_0)^{n}}\int_{\mathbb{B}^{n}}\frac{1 }{|x-y|^{n-2}}\, d\nu(y).
 \end{eqnarray*}
 Thus in order to prove the uniform convergence of $G_{h} [\psi](x)$ in $\overline{\mathbb{B}^{n}}(0,r_{0})$, we only need to prove that
 \[
F_{n-2}(x)=\int_{\mathbb{B}^{n}}\frac{1 }{|x-y|^{n-2}}\, d\nu(y)
\]
is uniformly convergent. In fact, we shall prove the following more general result.

\bcl\label{jl25-1cj} The integral $F_{k}(x)= \int_{\mathbb{B}^{n}}\frac{1 }{|x-y|^{k}}\, d\nu(y)$ is uniformly convergent w.r.t. $x$ in $\overline{\mathbb{B}^{n}}(0,r_{0})$, where $1\leq k\leq n-1$ and $0<r_{0}<1$. \ecl

Let $\delta_{2}=\frac{1-r_{0}}{3} $. Then
$$\mathbb{B}^{n}=\mathbb{B}^{n}(x,\delta_{2})\cup \big(\mathbb{B}^{n}\setminus \mathbb{B}^{n}(x,\delta_{2})\big)\;\;\mbox{and} \;\;\mathbb{B}^{n}(x,\delta_{2}) \subset \mathbb{B}^{n}.$$
Hence $$F_{k}(x)=F_{k,1}(x)+F_{k,2}(x),$$ where
$$F_{k,1}(x)=\int_{ \mathbb{B}^{n} \setminus \mathbb{B}^{n}(x,\delta_{2})}\frac{1 }{|x-y|^{k}}\, d\nu(y)\;\;\text{and}\;\;F_{k,2}(x)=\int_{\mathbb{B}^{n}(x,\delta_{2}) }\frac{1 }{|x-y|^{k}}\, d\nu(y).$$

\bscl\label{jl25-1c}  $F_{k,1}(x)$ and $F_{k,2}(x)$ are uniformly convergent w.r.t. $x$ in $\overline{\mathbb{B}^{n}}(0,r_{0})$, where $1\leq k\leq n-1$ and $0<r_{0}<1$. \escl

Since for all $y\in\mathbb{B}^{n} \setminus \mathbb{B}^{n}(x,\delta_{2})$,
$$\frac{1 }{|x-y|^{k}}\leq \frac{1 }{\delta_{2}^{k}},$$ by the Weierstrass test for uniform convergence, the uniform convergence of $F_{k,1}(x)$ in $\overline{\mathbb{B}^{n}}(0,r_{0})$ is obvious.

For any $0<\delta\leq \delta_2$, let $y=x+w$. Then it follows from \eqref{eq0.0cj.0} that
 $$\int_{\mathbb{B}^{n}(0,\delta)}\frac{1 }{|w|^{k}} \,d\nu(w)=\frac{n}{n-k}\delta^{n-k}\leq n\delta.$$
By definition, we easily know that $F_{k,2}(x)$ is uniformly convergent w.r.t. $x$ in $\overline{\mathbb{B}^{n}}(0,r_{0})$. Hence Subclaim \ref{jl25-1c} is proved.
 \medskip

Subclaim \ref{jl25-1c} implies the uniform convergence of $F_{k}(x)$ in $\overline{\mathbb{B}^{n}}(0,r_{0})$, and thus the proof of Claim \ref{jl25-1cj} is complete.\medskip

Let $k=n-2$. Then by Claim \ref{jl25-1cj}, we know that $G_{h} [\psi](x)$ is also uniformly convergent in $\overline{\mathbb{B}^{n}}(0,r_{0})$, and so the lemma is proved.  \epf

Now, we are going to prove the first main lemma for the proof of Theorem \ref{thm-2}.

\begin{lem}\label{lem-35.3}  Suppose $n\geq3$, $k\in\{1,\ldots,n\}$, $\psi\in C (\mathbb{B}^{n},\mathbb{R}^{n})$ and
$|\psi(x)|\leq M(1-|x|^{2}) $ in $\mathbb{B}^{n}$, where $M$ is a constant. Then
\begin{enumerate}
 \item for $0<r_{0}<1$, the unbounded integral $$I_{2,k}(x) =:\int_{\mathbb{B}^{n}} \left|\frac{\partial}{\partial x_{k}}G_{h}(x,y) \psi(y)\right| \,d\tau (y)$$ is uniformly convergent w.r.t. $x$ in $\overline{\mathbb{B}^{n}}(0,r_{0})$;
 \item for all $x\in   \mathbb{B}^{n} $, there exists a constant $\beta_{1}=\beta_{1}(n ,M)$ such that $$I_{2,k}(x)\leq \beta_{1}.$$
\end{enumerate}
\end{lem}

\bpf First, we easily see from \eqref{eq1.0.2}, \eqref{eq1.0.3}, \eqref{eq1.0.5} and \eqref{eq0.0.210} that for $x\neq y$,
\beq\label{eqcl4.2.2}
\frac{\partial}{\partial x_{k}}G_{ h}(x,y)
    &=& -\frac{(x_{k}-y_{k})(1-|x|^{2})^{n-1}(1-|y|^2)^{n-1}}{ n|x-y|^{n}[x,y]^{n}} \\ \nonumber
    &\;\;&-  \frac{x_{k}(1-|x|^{2})^{n-2}(1-|y|^2)^{n-1} }{ n|x-y|^{n-2}[x,y]^{n}  }.
 \eeq
 Then \eqref{wed-1} implies that for $x\in \mathbb{B}^{n}(0,r_{0})$, $[x,y]\geq 1-r_{0}$, and hence
 \begin{eqnarray*}
 \left|\frac{\partial}{\partial x_{k}}G_{ h}(x,y)\right|
    &\leq& \frac{ (1-|x|^{2})^{n-1}(1-|y|^2)^{n-1}}{ n|x-y|^{n-1}[x,y]^{n}} +\frac{ (1-|x|^{2})^{n-2}(1-|y|^2)^{n-1} }{ n|x-y|^{n-2}[x,y]^{n}  }\\
    &\leq& \frac{(1-|y|^2)^{n-1}}{ n(1-r_{0})^{n}} \left( \frac{1}{|x-y|^{n-1}}+ \frac{1}{|x-y|^{n-2}}\right) .\\\end{eqnarray*}
Thus the assumption ``$|\psi(x)|\leq M(1-|x|^2)$" implies that for all $x\in \mathbb{B}^{n}(0,r_{0})$,
 \begin{eqnarray*}
 I_{2,k}(x) &\leq&M\int_{\mathbb{B}^{n}}\left| \frac{\partial}{\partial x_{k}}G_{h}(x,y)\right|(1-|y|^{2})\,d\tau (y)\\
&\leq& \frac{M}{ n(1-r_{0})^{n}}\int_{\mathbb{B}^{n}}\left( \frac{1}{|x-y|^{n-1}}+ \frac{1}{|x-y|^{n-2}}\right) \,d\nu(y).\end{eqnarray*}
The uniform convergence of $I_{2,k}(x) $  w.r.t. $x$ in $\overline{\mathbb{B}^{n}}(0,r_{0})$ follows from Claim \ref{jl25-1cj}, and thus Lemma \ref{lem-35.3}(1) holds.

\medskip

Next, we prove Lemma \ref{lem-35.3}(2). It follows from \eqref{eqcl4.2.2} that
 \be\label{jl-00}\nonumber
 I_{2,k}(x)\leq \frac{1}{n}\big(I_{3,k}(x)+I_{4,k}(x)\big),\ee
 where
$$I_{3,k}(x)= \int_{\mathbb{B}^{n}}\frac{|x_{k}-y_{k}|(1-|x|^{2})^{n-1}}{ |x-y|^{n}[x,y]^{n}(1-|y|^2)}|\psi(y)|\,d\nu(y)$$
and
$$ I_{4,k}(x)= \int_{\mathbb{B}^{n}}\frac{|x_{k}|(1-|x|^{2})^{n-2} }{ |x-y|^{n-2}[x,y]^{n}(1-|y|^2) }  |\psi(y)|\,d\nu(y).$$

Next, we estimate $|I_{3,k}(x)|$ and $|I_{4,k}(x)|$, respectively.

 \bcl\label{jl24-3} For $x\in \mathbb{B}^{n}$,  $$|I_{3,k}(x)|\leq \frac{ n M  }{2} \mu_{2,3},$$
 where $\mu_{2,3}=\mu_{2}(n,\frac{1}{2},3,\frac{1}{2})$ is the same constant as in Lemma \ref{lem-4.5}. \ecl
Let $y=\varphi_{x}(w)$. Then the equalities \eqref{eq0.0cj.0}, \eqref{eq1.000.3}, \eqref{eq1.0.4} and the assumption
``$|\psi(x)|\leq M (1-|x|^{2})$" imply that
\begin{eqnarray*} \label{eq1.3.2}
\;\;\;|I_{3,k}(x)|&\leq& M \int_{\mathbb{B}^{n} }\frac{ (1-|x|^{2})^{n-1}|x-y|}{ |x-y|^{n}[x,y]^{n}}\,d\nu(y) \\
&=& M\int_{\mathbb{B}^{n} }\frac{ (1-|x|^{2})^{n-1}J_{\varphi_{x}}(w)}{ \left|x-\varphi_{x}(w)\right|^{n-1}[x,\varphi_{x}(w)]^{n}
}\,d\nu(w)\;\;\;\;\;\;\;\;\;\;\;\;\;\;\;\;\;\;\;\;\;\;\;\text{(by $y=\varphi_{x}(w)$)}\\
&=&  M \int_{\mathbb{B}^{n} }\frac{ d\nu(w)}{[x,w] \cdot|w|^{n-1}}\;\;\;\;\;\ \;\;\;\;\;\;\;\;\;\;\;\;\;\;\;\;\;\;\;\;\;\;\;\;\;\;\;\;\;\;\;\;\;\;\;\text{(by \eqref{eq1.000.3} and \eqref{eq1.0.4}))}\\
&=&n M \int_{0}^{1} \int_{\mathbb{S}^{n-1} }\frac{\,d\sigma(\xi)}{[x,\rho\xi]}\,d\rho.\;\;\;\;\;\;\;\;\;\;\;\;\;\;\;\;\;\;\;\;\;\;\;\;\;\;\;\;\;\;\;\;\;\;\;\;\;\;\;\;\;\;\;\;\;\;\;\;\;\;\;\;\;\;\;\;\;\text{(by \eqref{eq0.0cj.0})}
\end{eqnarray*}
Moreover, by \eqref{eq1.0.2}, we have
 $$\int_{\mathbb{S}^{n-1} }\frac{d\sigma(\xi)}{[x,\rho\xi] }=\int_{\mathbb{S}^{n-1} } \left(1+\rho^{2}|x|^{2}-2\rho|x|\left\langle\frac{x}{|x|},\xi\right\rangle\right)^{-\frac{1}{2}}\,d\sigma(\xi),$$
which, together with Theorem \Ref{Lem-D} and \eqref{eq2.2.4}, implies that
\begin{eqnarray*} \label{eq1.3.3}%\nonumber\;\;\;\;\;\;\;\;
\int_{\mathbb{S}^{n-1} }\frac{\,d\sigma(\xi)}{[x,\rho\xi] }
&=& \frac{\Gamma(\frac{n}{2})}{\Gamma(\frac{n-1}{2})\Gamma(\frac{1}{2})}\int_{-1}^{1} (1-s^{2})^{\frac{n-3}{2}}\left(1+\rho^{2}|x|^{2}
-2\rho|x|s\right)^{-\frac{1}{2}}ds\\
 &=&F\left(\frac{ 1}{2},\frac{3-n}{2};\frac{n}{2};\rho^{2}|x|^{2}\right).
\end{eqnarray*}
Hence Lemma \ref{lem-4.2} leads to
\begin{eqnarray*}
  |I_{3,k}(x)| &\leq & nM   \int_{0 }^{1}  F\left(\frac{ 1}{2},\frac{3-n}{2};\frac{n}{2};\rho^{2}|x|^{2}\right)\,d\rho\leq \frac{ n M  }{2}\mu_{2,3},
 \end{eqnarray*}
as required, where $\mu_{2,3}=\mu_{2}(n,\frac{1}{2},3,\frac{1}{2})$.

\bcl\label{jl25-1}  For $x\in \mathbb{B}^{n}$, we have  $$|I_{4,k}(x)|\leq M\mu_5, $$ where $\mu_5=\mu_5(n)$ is the same constant as in Lemma \ref{lem1-100}. \ecl
Obviously, the assumption ``$|\psi(x )|\leq M (1- |x|^2 )$" implies that
 $$|I_{4,k}(x)|\leq  M\int_{\mathbb{B}^{n}} \frac{ (1-|x|^{2})^{n-2} }{ |x-y|^{n-2}[x,y]^{n} }\,d\nu(y).$$
Let $y=\varphi_{x}(w)$. By \eqref{eq1.000.3}, \eqref{eq1.0.4} and Lemma \ref{lem1-100}, we get
\begin{eqnarray*}
\;\;\;\; |I_{4,k}|
&\leq&   M\int_{\mathbb{B}^{n}} \frac{ (1-|x|^{2})^{n-2}J_{\varphi_{x}}(w)}{ |x-\varphi_{x}(w)|^{n-2}[x,\varphi_{x}(w)]^{n}}\,d\nu(w)\;\;\;\;\;\;\; \;\;\;\text{(substituting $y=\varphi_{x}( w)$)}
\\%\nonumber
&=&M\int_{\mathbb{B}^{n} }\frac{ d\nu(w)}{|w|^{n-2}[x,w]^{2}}\leq M\mu_5(n).
 \end{eqnarray*}

By taking $\beta_{1} =\frac{ M  }{2}\mu_{2,3} +\frac{M}{n}\mu_5 $, we see that Lemma \ref{lem-35.3}(2) holds, and so the proof of the lemma is finished.
\epf

\begin{lem}\label{lem-34.3}  Suppose $n\geq3$, $\psi\in C (\mathbb{B}^{n},\mathbb{R}^{n})$ and
$|\psi(x)|\leq M(1-|x|^{2}) $ in $\mathbb{B}^{n}$, where $M$ is a constant. Then for all $0<r_0<1$
and $k\in \{1,2,\ldots,n\}$, $$G_{h}[\psi](x)\;\;\mbox{and}\;\; \int_{\mathbb{B}^{n}} \frac{\partial}{\partial x_{k}} G_{h}(x,y)\psi(y)\,d\tau(y)$$ are continuous in $\overline{\mathbb{B}^{n}}(0,r_{0})$, respectively.
\end{lem}
\bpf In order to check the continuity of $G_{h}[\psi](x)$ in $\overline{\mathbb{B}^{n}}(0,r_{0})$, we only need to prove that
$G_{h}[\psi](x)$ is continuous at every fixed point $x_0 \in \overline{\mathbb{B}^{n}}(0,r_{0})$.
Assume that $x_{0}\in \overline{\mathbb{B}^{n}}(0,r_{0})$ and $x_{0}+\Delta x\in \overline{\mathbb{B}^{n}}(0,r_{0})$.

By Lemma \ref{lem-4.23}, we see that $ G_{h} [\psi](x)$ is uniformly convergent in $\overline{\mathbb{B}^{n}}(0,r_{0})$.
Then for any $\varepsilon_{1}>0$, there exist constants
$\iota_{1}= \iota_{1}(\varepsilon_{1})\rightarrow 1^{-}$ and $\iota_{2}=\iota_{2}(\varepsilon_{1})\rightarrow 0^{+}$ such that for any $x\in \overline{\mathbb{B}^{n}}(0,r_{0})$,
$$\mathbb{B}^{n}(x,\iota_{2} ) \subset \mathbb{B}^{n}(0,\iota_{1}),$$
 $$\left|\int_{\mathbb{B}^{n}\backslash \mathbb{B}^{n}(0,\iota_{1} )}  G_{h}(x,y)\psi(y) \,d\tau(y)\right|<\varepsilon_{1} \;\;
 \text{and}
\;\; \left|\int_{\mathbb{B}^{n}(x,\iota_{2} )}  G_{h}(x,y)\psi(y) \,d\tau(y)\right|<\varepsilon_{1}.$$
 Then
\beq \label{cjeq1.3.200}
&\;\;& \left|G_{h}[\psi](x_{0}+\Delta x)-G_{h}[\psi](x_{0})\right|\\ \nonumber
&=& \left|\int_{\overline{\mathbb{B}^{n}}(0,\iota_{1})\backslash \mathbb{B}^{n}(x_{0},\iota_{2})}  \big[ G_{h}(x_{0}+\Delta x,y)-G_{h}(x_{0},y)\big]\psi(y) \,d\tau(y) \right.\\ \nonumber
&\;\;&+\int_{\mathbb{B}^{n}(x_{0},\iota_{2})} \big[ G_{h}(x_{0}+\Delta x,y)-G_{h}( x_{0},y)\big]\psi(y) \,d\tau(y)\\ \nonumber
&\;\;&\left.+\int_{\mathbb{B}^{n}\backslash  \mathbb{B}^{n}(0,\iota_{1})}\big[  G_{h}(x_{0}+\Delta x,y)- G_{h}(x_{0} ,y)\big]\psi(y) \,d\tau(y)\right|\\ \nonumber
&\leq &\left| \int_{\overline{\mathbb{B}^{n}}(0,\iota_{1})\backslash \mathbb{B}^{n}(x_{0},\iota_{2})}\big[ G_{h}(x_{0}+\Delta x,y)-G_{h}(x_{0},y)\big]\frac{\psi(y)}{(1-|y|^{2})^{n}}\,d\nu(y) \right| +4\varepsilon_{1}.
\eeq

By \eqref{eq0.0.210}, it is easy to see that the map $(x,y)\rightarrow G_{h}(x,y)\frac{\psi(y)}{(1-|y|^{2})^{n}}$ is continuous (also uniformly continuous) on $ \overline{\mathbb{B}^{n}}(x_{0}, \frac{1}{2}\iota_{2})\times  \big(\overline{\mathbb{B}^{n}}(0,\iota_{1})\backslash \mathbb{B}^{n}(x_{0},\iota_{2})\big)$. Therefore, there exists $\iota'=\iota'(\varepsilon_{1})<\frac{1}{2}\iota_{2}$
such that for all $|\Delta x|<\iota'$
 and for all $y\in \overline{\mathbb{B}^{n}}(0,\iota_{1})\backslash \mathbb{B}^{n}(x_{0},\iota_{2}) $,
\be\label{eqcjl4.2.5} \left|\big(G_{h}(x_{0}+\Delta x,y)-G_{h}(x_{0},y)\big)\frac{\psi(y)}{(1-|y|^{2})^{n}}\right|<\varepsilon_{1} .\ee
Thus it follows from \eqref{cjeq1.3.200} and \eqref{eqcjl4.2.5} that
 $$ \left|G_{h}[\psi](x_{0}+\Delta x)-G_{h}[\psi](x_{0})\right|\leq5\varepsilon_{1},$$
which means that $G_{h}[\psi]$ is continuous at $x_0$. Hence the arbitrariness of $x_{0}$ shows that $G_{h}[\psi](x)$ is continuous in $\overline{\mathbb{B}^{n}}(0,r_{0})$.

By applying Lemma \ref{lem-35.3}, the continuity of $$\int_{\mathbb{B}^{n}} \frac{\partial}{\partial x_{k}} G_{h}(x,y)\psi(y)\,d\tau(y)$$ in $\overline{\mathbb{B}^{n}}(0,r_{0})$ can be proved in a similar way as above, where $k\in \{1,\;\ldots, \;n\}$. So the proof of this lemma is complete.
\epf

The following property is the second main lemma for the proof of Theorem \ref{thm-2}.

\begin{lem}\label{lem-4.3}  Suppose $n\geq3$, $\psi\in C (\mathbb{B}^{n},\mathbb{R}^{n})$ and
$|\psi(x)|\leq M(1-|x|^{2}) $ in $\mathbb{B}^{n}$, where $M$ is a constant. Then for all $x\in \mathbb{B}^{n}$
and $k\in \{1,2,\ldots,n\}$,
%\be\label{eq4.2.0}
$$\frac{\partial \Psi}{\partial x_{k}}(x)
 =\int_{\mathbb{B}^{n}}  \frac{\partial}{\partial x_{k}}G_{h}(x,y)\psi(y) \,d\tau(y).$$
\end{lem}

 \bpf For all $x\in \mathbb{B}^{n}$, by Lemma \ref{lem-35.3}, we see that $$\int_{0}^{x_{k}}\int_{\mathbb{B}^{n}} \left|\frac{\partial}{\partial x_{k}} G_{h}(x,y)\frac{\psi(y)}{(1-|y|^2)^{n-1}}\right|\,d\nu(y)dx_{k}\leq \beta_{1},$$ where $k\in \{1,\;\ldots, \;n\}$.
It follows from Fubini's theorem \cite[p. 165]{rudin3} that
$$\int_{0}^{x_{k}}\int_{\mathbb{B}^{n}} \frac{\partial}{\partial x_{k}} G_{h}(x,y)\psi(y)\,d\tau(y)dx_{k}
  = \int_{\mathbb{B}^{n}} \int_{0}^{x_{k}}\frac{\partial}{\partial x_{k}}G_{h}(x,y)\frac{\psi(y)}{(1-|y|^{2})^{n}} \,dx_{k}\,d\nu(y),
$$
which means
 \begin{eqnarray*}
&\;\;&\int_{0}^{x_{k}}\int_{\mathbb{B}^{n}} \frac{\partial}{\partial x_{k}} G_{h}(x,y)\psi(y)\,d\tau(y)dx_{k}\\ \nonumber
&=& \int_{\mathbb{B}^{n}} \frac{G_{h}(x,y) }{(1-|y|^{2})^{n}} \psi(y) \,d\nu(y)-\int_{\mathbb{B}^{n}} \frac{ G_{h}(x_{k,0},y)}{(1-|y|^{2})^{n}} \psi(y) \,d\nu(y),
 \end{eqnarray*}
where $x_{k,0}=(x_{1},\ldots,x_{k-1},0,x_{k+1},\ldots,x_{n})$.
 Since $\int_{\mathbb{B}^{n}} \frac{\partial}{\partial x_{k}} G_{h}(x,y)\psi(y)\,d\tau(y)$ is continuous in $\overline{\mathbb{B}^{n}}(0,r_{0})$, by differentiating w.r.t. $x_{k}$, we get
 $$\int_{\mathbb{B}^{n}} \frac{\partial}{\partial x_{k}} G_{h}(x,y)\psi(y)\,d\tau(y)= \frac{\partial}{\partial x_{k}} \int_{\mathbb{B}^{n}} \frac{G_{h}(x,y) }{(1-|y|^{2})^{n}} \psi(y) \,d\nu(y).$$
 Hence the proof of this lemma is finished. \epf

\subsection*{Proof of Theorem \ref{thm-2}}
Lemmas \ref{lem-35.3}(2), \ref{lem-4.3} and Cauchy-Schwarz inequality imply that
%\be\label{eq100.300.300}
$$\|D \Psi(x)\|=\sup_{\xi\in \mathbb{S}^{n-1}} |D \Psi(x)\xi |\leq  \left(\sum_{k=1}^{n}\left| \frac{\partial \Psi}{\partial x_{k}}(x) \right|^{2}  \right)^{\frac{1}{2}}\leq \sqrt{n}\beta_{1}=:\beta_{0},$$
where $\beta_{1}=\beta_{1}(n,M)$ is the same constant as in Lemma \ref{lem-35.3}(2). For any $x,y\in \mathbb{B}^{n}$, let $\gamma_{[x,y]}$ denote the segment between $x$ and $y$.
By Lemmas \ref{lem-34.3} and \ref{lem-4.3}, we know that $\Psi\in C^{1}( \mathbb{B}^{n})$, and hence $\Psi$ is differentiable.
Then the mean-value theorem of differentials leads to $$\left|\Psi(x)-\Psi(y)\right|\leq   \beta_{0} |x-y|.$$
The proof of Theorem \ref{thm-2} is finished.\qed
\medskip

Based on  Theorems \ref{thm-11} and \ref{thm-2}, we are going to prove Theorem \ref{thm-4}.

\subsection*{Proof of Theorem \ref{thm-4}}
In this subsection, we always regard a point $x=(x_1, \ldots, x_n) $ in $\mathbb{R}^n$ as a column vector, for purposes of computing matrix products (which have been denoted by $\times$).

For any $x,$ $y\in \mathbb{B}^{n}$, by letting $\omega(t)=Lt$ in Lemma \ref{jl23-3}, we obtain that for $x\in \mathbb{B}^{n}$,
%\be\label{sun-1}
$$\|D \Phi(x) \|\leq L\alpha_{0},$$
where $L$ is the same constant as in Theorem \ref{thm-4}.
It follows from the mean-value theorem of differentials that
$$|\Phi(x)-\Phi(y)|\leq L \alpha_{0}  |x-y|,$$
and so Theorem \ref{thm-2} gives
$$
 |u(x)-u(y)|\leq  |\Phi(x)-\Phi(y)|+ |\Psi(x)-\Psi(y)|
\leq(L \alpha_{0}+\beta_{0} ) |x-y|.
$$

Let $N=L \alpha_{0}+\beta_{0}$. Then the proof of Theorem \ref{thm-4} is complete.\qed
\medskip

\section{Examples}\label{example-sec}

In this section, we will first construct an example to show that the requirement $n\geq 3$ in Theorem \ref{thm-4} is necessary.

\beg\label{ex-1}
Let $w_{0}(re^{i\theta})=\sum_{k=1}^{\infty}\frac{ r^{k}}{k^{2}}\cos (k\theta)-\frac{M}{4}(1-r^2)$ in $\mathbb{D}$, where $M$ is a non-negative constant. Then,
\begin{enumerate}
\item
$w_{0}\in C^2(\mathbb{D}, \mathbb{R})\cap C(\overline{\mathbb{D}}, \mathbb{R})$ and $\Delta_{h} w_{0}=M(1-|z|^{2})^{2}$;
\item
$w_{0}$ is not Lipschitz continuous in $\mathbb{D}$;
\item$w_{0}|_{\mathbb{S}^{1}}$ is Lipschitz continuous in $\mathbb{S}^{1}$.
\end{enumerate}
\eeg
\bpf
To prove that the function $w_0$ has the desired properties, let $$f(z)=\sum_{k=1}^{\infty}\frac{z^{k} }{k^{2}}.$$
Then for $z\in \mathbb{D}$,
$$w_{0}(z)=\text{Re}f(z)-\frac{M}{4}(1-|z|^2).$$
Obviously, $\text{Re}f=P[\phi_{0}] $ is harmonic in $\mathbb{D}$,
where $\phi_{0}(e^{i\theta})=\sum_{k=1}^{\infty}\frac{1}{k^{2}}\cos (k\theta)$ is continuous in $\mathbb{S}^{1}$,
and thus  $\text{Re}f\in C^2(\mathbb{D}, \IR)\cap (\overline{\mathbb{D}}, \IR)$.  By elementary computations, we see that $\Delta w_{0}=M$.
Hence the first assertion in the example holds.

Since $$\frac{\partial }{\partial z}w_{0}(z)=\frac{1}{2}\sum_{k=1}^{\infty}\frac{z^{k-1}}{k}+\frac{M}{4}\overline{z}=-\frac{\log(1-z)}{2z}+\frac{M}{4}\overline{z}$$
and
$$\frac{\partial }{\partial \overline{z}}w_{0}(z)=\frac{1}{2}\sum_{k=1}^{\infty}\frac{\overline{z}^{k-1}}{k}+\frac{M}{4}z=-\frac{\log(1-\overline{z})}{2\overline{z}}+\frac{M}{4}z,$$
we easily see that
 $$\|D w_{0}(z)\|=  \left|\frac{\partial }{\partial z}w_{0}(z)\right|+\left|\frac{\partial }{\partial \overline{z}}w_{0}(z)\right| $$ is unbounded in $\mathbb{D}$.
\bcl\label{sat-1}
The function $w_{0}$ is Lipschitz continuous if and only if $\|Dw_{0}\|$ is bounded.
\ecl
For the proof, we let $\partial_{\theta}w_{0}(z)$ denote the directional derivative of $w_{0}$.
  If $w_{0}$ is Lipschitz continuous with Lipschitz constant $L_{1}$, then
$$\big | \partial_{\theta}w_{0}(z) \big|=\left|\lim_{r\to 0}\frac{w_{0}(z+re^{i \theta})-w_{0}(z)}{r } \right|= \lim_{r\to0}\frac{\left|w_{0}(z+re^{i \theta})-w_{0}(z)\right|}{r } \leq L_{1}.$$ Hence it follows from the obvious fact
%{\color{red}  (Since $w_{0}(z+re^{i\theta})-w_{0}(z)=\frac{\partial}{\partial z}w_{0}(z)re^{i\theta}+\frac{\partial}{\partial
%\overline{z}}w_{0}(z)\overline{re^{i\theta}}+o(r^{2})$)}
$\|Dw_{0}(z)\|=\max_{\theta}\big | \partial_{\theta}w_{0}(z) \big|$ that $$\|Dw_{0}(z)\|\leq L_{1}.$$

On the other hand, if $\|Dw_{0}(z)\|\leq L_{1}$, then the mean-value theorem of differentials leads to
$$|w_{0}(z_{1})-w_{0}(z_{2})|\leq L_{1}|z_{1}-z_{2}|.$$
Hence the claim is true.\medskip

Since we have proved that $\|D w_{0}(z)\|$ is unbounded, we see from Claim \ref{sat-1} that $w_{0}$ is not Lipschitz continuous in $\mathbb{D}$, which shows that the second assertion in the example holds too. The third assertion follows from \cite[p. 317]{AKM2008} as the construction of $w_0$ in $\mathbb{S}^{1}$ coincides with the one in \cite{AKM2008}.
\epf

\medskip

Now, we construct an example to show that there exists a hyperbolic harmonic mapping $u$, which satisfies the conditions of Theorem \ref{thm-11} but $u$ is not $K$-quasiregular.
Recall that Example 1 
of \cite{AKM2008} shows that Lipschitz continuity of $\phi:\mathbb{S}^{n-1}\rightarrow \mathbb{R}^{n}$ does not imply Lipschitz continuity of its harmonic extension $u=P[\phi]:\mathbb{B}^{n}\rightarrow \mathbb{R}^{n}$, without the assumption that $u$ is quasiregular.

\beg\label{ex-2}
For $n=4$, let $u(x)=(1-\frac{1}{3}|x|^{2})x$ in $\mathbb{B}^{4}$. Then
\begin{enumerate}
\item $u \in C^2(\mathbb{B}^{4}, \mathbb{R}^{4})$ $\cap$ $C(\overline{\mathbb{B}^{4}}, \mathbb{R}^{4})$ and $\Delta_{h}u =0$ in $\mathbb{B}^{4}$;
\item $u=P_{h}[\phi]$, where $\phi(\xi)=\frac{2}{3}\xi$ and $\xi \in \mathbb{S}^{3}$;
\item $u$ is not $K$-quasiregular for any $K\geq1$;
\item $u$ is Lipschitz continuous in $\mathbb{B}^{4}$;
\item $\phi$ is Lipschitz continuous in $\mathbb{S}^{3}$.

\end{enumerate}
\eeg
\bpf
Elementary calculations and Theorem \Ref{Thm-B} show that $\Delta_{h}u =0$ and $u$ satisfies the conditions (1) and (2).

For any $x=(x_{1},\cdots, x_{4})\in \mathbb{B}^{4}$, it follows from \eqref{ja} that
$$Du(x)=\left(
             \begin{array}{cccc}
               1-\frac{1}{3}|x|^{2}-\frac{2}{3}x_{1}^{2} & -\frac{2}{3}x_{1}x_{2} & -\frac{2}{3}x_{1}x_{3} &-\frac{2}{3}x_{1}x_{4} \\
               \; & \; &\;&\;\\
               -\frac{2}{3}x_{1}x_{2} & 1-\frac{1}{3}|x|^{2}-\frac{2}{3}x_{2}^{2} & -\frac{2}{3}x_{2}x_{3} & -\frac{2}{3}x_{2}x_{4} \\
               \; & \; &\;&\;\\
               -\frac{2}{3}x_{1}x_{3} & -\frac{2}{3}x_{2}x_{3} & 1-\frac{1}{3}|x|^{2}-\frac{2}{3}x_{3}^{2} & -\frac{2}{3}x_{3}x_{4} \\
               \; & \; &\;&\;\\
              -\frac{2}{3}x_{1}x_{4} & -\frac{2}{3}x_{2}x_{4} & -\frac{2}{3}x_{3}x_{4} & 1-\frac{1}{3}|x|^{2}-\frac{2}{3}x_{4}^{2} \\
             \end{array}
           \right)
.$$
By calculations, we see that the eigenvalues of $Du(x)^T\times Du(x)$ are
$$\lambda_{1}^{2}= (1- |x|^2)^2,\;\;\lambda_{2}^{2}=\lambda_{3}^{2}=\lambda_{4}^{2}=\frac{1}{9} (3- |x|^2)^2.$$
Thus, $$J_{u}(x)=\frac{1}{27}(3- |x|^2)^3(1- |x|^2)\;\;\text{ and}\;\;||Du(x)||=1-\frac{1}{3}|x|^{2},$$
which implies that $u$ is not $K$-quasiregular in $ \mathbb{B}^{4}$, and the statement (3) holds.

Further, for any $x,y \in \mathbb{B}^{4}$ and $\xi,\eta \in \mathbb{S}^{3}$, we have
 \begin{eqnarray*}
|u(x)-u(y)|&\leq& |x-y|+\frac{1}{3}\big|  |x|^{2}x-|y|^{2}y\big|\\ \nonumber
&\leq& |x-y|+\frac{1}{3} |x|^{2}\cdot|x-y|+\frac{1}{3}( |x|^{2}- |y|^{2}) |y|\\ \nonumber
&\leq& 2|x-y|,
 \end{eqnarray*}
and $ |\phi(\xi)-\phi(\eta)|=\frac{2}{3}|\xi-\eta| $.
Thus, the statements (4) and (5) hold.
  \epf

\medskip

\section*{Acknowledgements}
 We thank Professor Manfred Stoll for useful suggestions to this paper.

\end{document}